\documentclass[a4paper,10pt]{article}

\usepackage{amsmath,amsfonts,amssymb,amsthm}
\usepackage{mathtools}
\usepackage{stmaryrd}
\usepackage{latexsym}
\usepackage{cite}
\usepackage{siunitx}
\usepackage{wasysym}
\usepackage{url}
\usepackage{xspace}
\usepackage{graphicx}

\usepackage[small,bf]{caption}
\usepackage{multirow}
\usepackage{float}
\usepackage{subfig}

\usepackage{listings}
\usepackage{xcolor}
\usepackage{longtable}
\usepackage{listings}
\lstdefinelanguage{pseudo}{
frame=BT,
mathescape,
morekeywords={1.,2.,3.,4.,5.,6.,7.,8.,9.,10.,11.,12.,(a),(b),(c),(d)},
basicstyle=\normalsize \ttfamily \color{black},
showspaces=false, 
showstringspaces=false,        
showtabs=false,
keywordstyle=\bfseries,
}

\def \HDiv {H(\operatorname{div})}
\def \Hdiv0 {H_0(\operatorname{div})}

\newcommand{\Corg}{C^{\textrm{org}}}
\newcommand{\Cdx}{C^{\textrm{dx}}}
\newcommand{\Fout}{F^{\textrm{out}}}
\newcommand{\Topt}{T_\textrm{opt}}
\newcommand{\wMax}{w_\textrm{max}}

\newcommand{\Prod}{\Theta_g}
\newcommand{\Nprod}{N_g}
\newcommand{\Jout}{J_\textrm{out}}
\newcommand{\Inj}{\Theta_w}
\newcommand{\Ninj}{N_w}
\newcommand{\Pvap}{P^{\textrm{vp}}}
\newcommand{\Fcond}{F^{\textrm{cond}}}
\newcommand{\Fevap}{F^{\textrm{evap}}}
\newcommand{\GammaTop}{\Gamma_\textrm{t}}
\newcommand{\GammaBase}{\Gamma_\textrm{b}}
\newcommand{\GammaLat}{\Gamma_\textrm{l}}
\newcommand{\Tm}{T_{\textrm{m}}}
\newcommand{\Tg}{T_{\textrm{g}}}
\newcommand{\Jin}{J_\textrm{in}}
\newcommand{\Fin}{F^{\textrm{in}}}
\newcommand{\Lgls}{\mathcal{L}^{\textrm{GLS}}}
\newcommand{\Lsupg}{\mathcal{L}^{\textrm{SUPG}}}
\newcommand{\Line}{\mathfrak{L}}
\newcommand{\Sfin}{S^{\textrm{fin}}}
\newcommand{\Feel}{\textsc{Feel++}\xspace} 
\newcommand{\Gmsh}{\textsc{Gmsh}\xspace}
\numberwithin{equation}{section}
\theoremstyle{plain}
\newtheorem{theorem}{Theorem}[section]
\theoremstyle{remark}
\newtheorem{rmrk}[theorem]{Remark}
\usepackage[hang,flushmargin]{footmisc}
\providecommand{\keywords}[1]{{\small{\textbf{Keywords:}} #1}}

\title{Mathematical modeling and numerical simulation of a bioreactor landfill using Feel++}

\author{
\renewcommand{\thefootnote}{\alph{footnote}}
G. Doll\'e \footnotemark[1] , O. Duran \footnotemark[2] , N. Feyeux 
\footnotemark[3] , E. Fr\'enod \footnotemark[4] , M. Giacomini 
\footnotemark[5]\textsuperscript{ \ ,}\footnotemark[6] \ and C. Prud'homme 
\footnotemark[1]
}

\date{}

\begin{document}

\renewcommand{\thefootnote}{\alph{footnote}}

\maketitle

\footnotetext[1]{Universit\'e de Strasbourg, IRMA UMR 7501, 7 rue Ren\'e 
Descartes, 67084 Strasbourg, France.}
\footnotetext[2]{State University of Campinas, SP, Brazil.}
\footnotetext[3]{MOISE team, INRIA Grenoble Rh\^one-Alpes. Universit\'e de 
Grenoble, Laboratoire Jean Kuntzmann, UMR 5224, Grenoble, France.}
\footnotetext[4]{Universit\'e de Bretagne-Sud,  UMR 6205, LMBA, F-56000 Vannes, 
France.}
\footnotetext[5]{CMAP, Inria, Ecole polytechnique, CNRS, Universit\'e Paris-Saclay, 91128 Palaiseau, France.}
\footnotetext[6]{DRI Institut Polytechnique des Sciences Avanc\'ees, 63 Boulevard de Brandebourg, 94200 Ivry-sur-Seine, France.\vspace{5pt} }

\footnotetext{
\textit{
This work has been supported by the LMBA Universit\'e de Bretagne-Sud, 
the project PEPS Amies VirtualBioReactor and the private funding of 
See-d and Entreprise Charier.  \\
The project is hosted on the facilities at CEMOSIS whose support is kindly
acknowledged. \\
M. Giacomini is member of the DeFI team at Inria Saclay \^Ile-de-France.
}
\vspace{5pt} }

\footnotetext{\texttt{\textit{e-mail:}
dolle@math.unistra.fr;
omar@dep.fem.unicamp.br;
nelson.feyeux@imag.fr;
emmanuel.frenod@univ-ubs.fr; 
matteo.giacomini@polytechnique.edu;
prudhomme@unistra.fr.
}
}

\renewcommand{\thefootnote}{\arabic{footnote}}

\begin{abstract}
In this paper, we propose a mathematical model to describe the functioning of a
bioreactor landfill, that is a waste management facility in which biodegradable
waste is used to generate methane.
The simulation of a bioreactor landfill is a very complex multiphysics problem
in which bacteria catalyze a chemical reaction that starting from
organic carbon leads to the production of methane, carbon dioxide and water.
The resulting model features a heat equation coupled with a non-linear reaction
equation describing the chemical phenomena under analysis and several advection
and advection-diffusion equations modeling multiphase flows inside a porous
environment representing the biodegradable waste.
A framework for the approximation of the model is implemented using \Feel, a
C++ open-source library to solve Partial Differential Equations.
Some heuristic considerations on the quantitative values of the parameters in
the model are discussed and preliminary numerical simulations are presented.
\end{abstract}

\maketitle
%
\section{Introduction}

Waste management and energy generation are two key issues in nowadays
societies. A major research field arising in recent years focuses on combining
the two aforementioned topics by developing new techniques to handle waste
and to use it to produce energy.
A very active field of investigation focuses on bioreactor landfills which are
facilities for the treatment of biodegradable waste. The waste is accumulated
in a humid environment and its degradation is catalyzed by bacteria. The
main process taking place in a bioreactor landfill is the methane generation
starting from the consumption of organic carbon due to waste decomposition.
Several by-products appear during this reaction, including carbon dioxide and
leachate, that is a liquid suspension containing particles of the waste
material through which water flows.

Several works in the literature have focused on the study of bioreactor
landfills but to the best of our knowledge none of them tackles the
global multiphysics problem.
On the one hand, \cite{PMC98, Reinhart1996, Martin2001} present
mathematical approaches to the problem but the authors deal with a single
aspect of the phenomenon under analysis focusing either on microbiota activity
and leachates recirculation or on gas dynamic.
On the other hand, this topic has been of great interest in the engineering
community \cite{TOUGH06, AIC:AIC690490821, Kindlein2006} and several studies
using both numerical and experimental
approaches are available in the literature.
We refer the interested reader to the review paper \cite{Agostini01102012} on
this subject.

In this work, we tackle the problem of providing a mathematical model for the
full multiphysics problem of methane generation inside a bioreactor landfill.
Main goal is the development of a reliable model to simulate the long-time
behavior of these facilities in order to be able to perform forecasts and
process optimization \cite{Liu2006}.
This paper represents a preliminary study of the problem starting from the
physics of the phenomena under analysis and provides a first set of equations
to describe the methane generation inside a bioreactor landfill.
In a more general framework, we aim to develop a model sufficiently accurate
to be applied to an industrial context limiting at the same time the required
computational cost. 
Thus, a key aspect of this work focused on the identification 
of the most important features of the functioning of a bioreactor landfill 
in order to derive the simplest model possible to provide an accurate 
description of the aforementioned methanogenic phenomenon.
The proposed model has been implemented using \Feel and
the resulting tool to numerically simulate the dynamic of a bioreactor landfill
has been named SiViBiR++ which stands for \emph{Simulation of a Virtual
BioReactor using \Feel}.

The rest of the paper is organized as follows. After a brief description of the
physical and chemical phenomena taking place inside these waste management
facilities (Section \ref{ref:physical_model}), in section
\ref{ref:math_model} we present the fully coupled mathematical model of a
bioreactor landfill. Section \ref{ref:numerical_model} provides details on the
numerical strategy used to discretize the discussed model. Eventually, in
section \ref{ref:simulations} preliminary numerical tests are presented and
section \ref{ref:conclusion} summarizes the results and highlights some future
perspectives. In appendix \ref{sec:parameters}, we provide a table with the 
known and unknown parameters featuring our model.

\section{What is a bioreactor landfill?}
\label{ref:physical_model}

As previously stated, a bioreactor landfill is a facility for the treatment of
biodegradable waste which is used to generate methane, electricity and hot water.
Immediately after being deposed inside a bioreactor, organic waste begins
to experience degradation through chemical reactions. During the first
phase, degradation takes place via aerobic metabolic pathways, that is a series
of concatenated biochemical reactions which occur within a cell in presence of
oxygen and may be accelerated by the action of some enzymes. Thus bacteria
begin to grow and metabolize the biodegradable material and complex organic
structures are converted to simpler soluble molecules by hydrolysis.

The aerobic degradation is usually short because of the high demand of oxygen
which may not be fulfilled in bioreactor landfills. Moreover, as more material
is added to the landfill, the layers of waste tend to be compacted and the
upper strata begin to block the flow of oxygen towards the lower parts
of the bioreactor. Within this context, the dominant reactions inside the
facility become anaerobic. Once the oxygen is exhausted, the bacteria begin to
break the previously generated molecules down to organic acids which are
readily soluble in water and the chemical reactions involved in the metabolism
provide energy for the growth of population of microbiota.

After the first year of life of the facility, the anaerobic conditions and the
soluble organic acids create an environment where the methanogenic bacteria can
proliferate \cite{WalshKinman}. These bacteria become the major actors inside
the landfill by using the end products from the first stage of degradation to
drive the methane fermentation and convert them into methane and carbon
dioxide. Eventually, the chemical reactions responsible for the generation of
these gases gradually decrease until the material inside the landfill is inert
(approximately after 40 years). \\
In this work, we consider the second phase of the degradation process, that is
the methane fermentation during the anaerobic stage starting after the first
year of life of the bioreactor.

\subsection{Structure of a bioreactor landfill}

\begin{figure}[tbp]
\centering
\subfloat[3D render of a bioreactor landfill.]
{
\includegraphics[width=0.36\columnwidth]{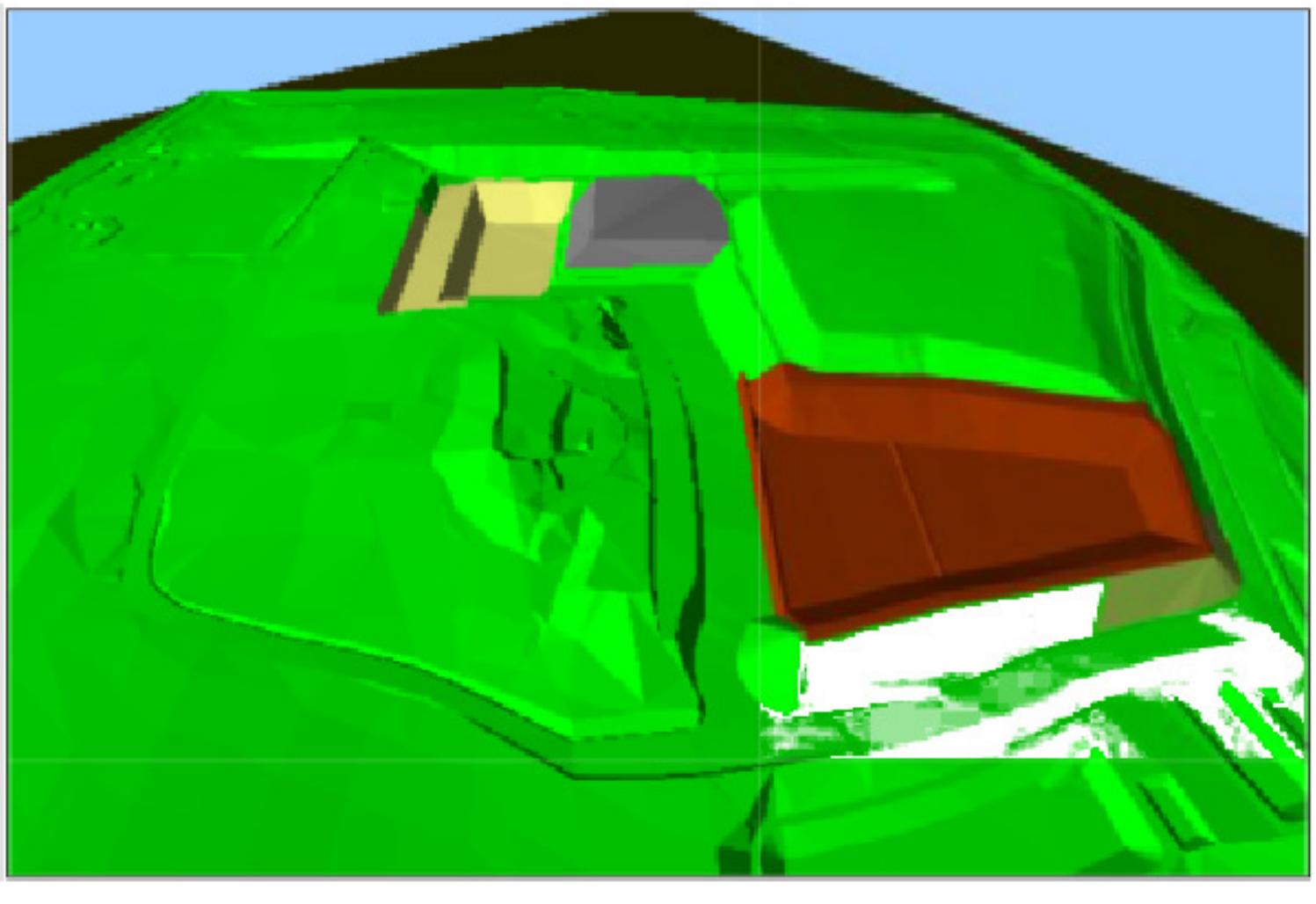}
\label{fig:droues}
}
\hfill
\subfloat[Scheme of the structure of an alveolus.]
{
\includegraphics[width=0.6\columnwidth]{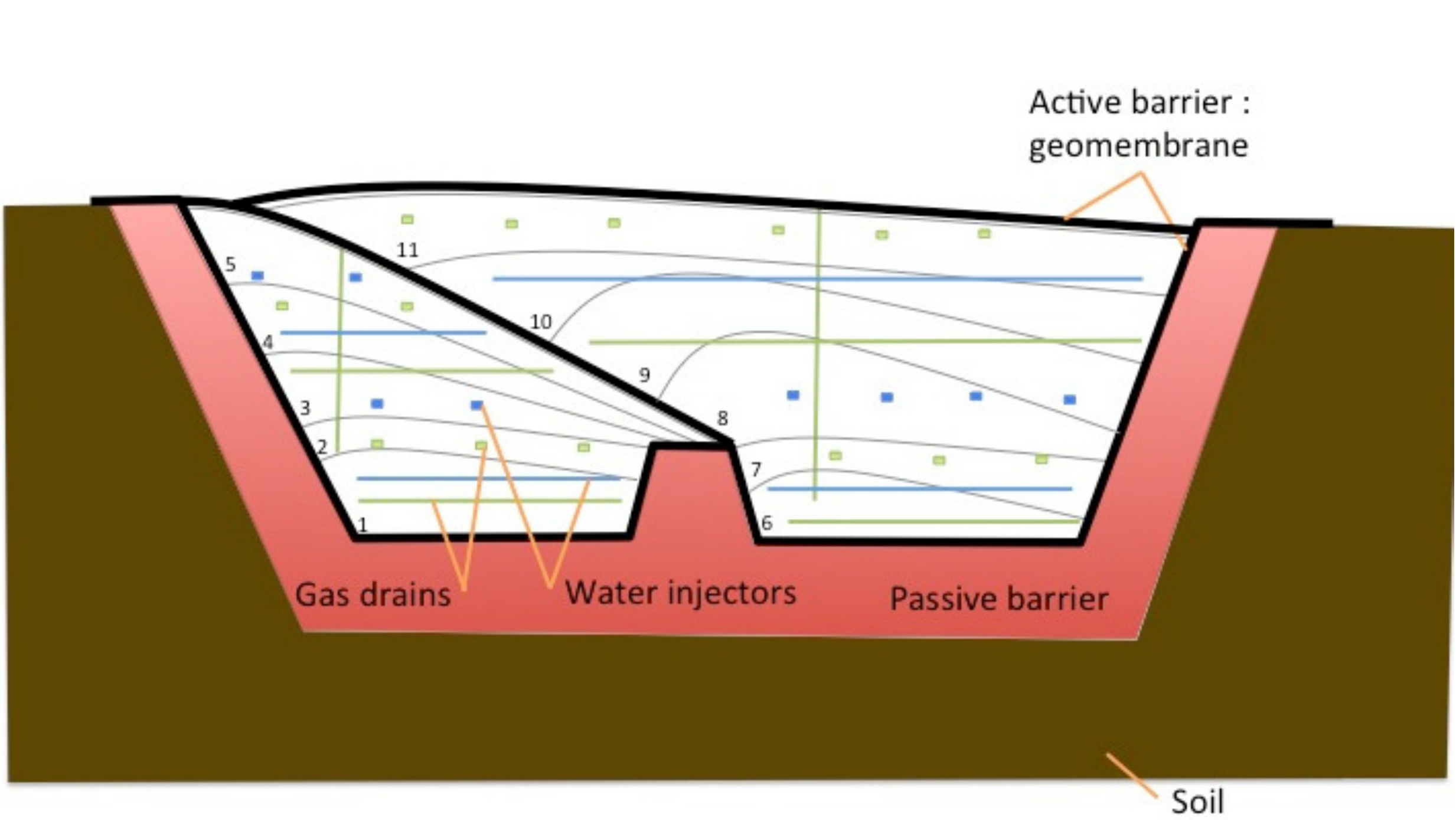}
\label{fig:alveolus_scheme}
}
\caption{Structure of a bioreactor landfill and its composing alveoli. Image
courtesy of Entreprise Charier \protect\url{http://www.charier.fr}.}
\label{fig:schemes}
\end{figure}

A bioreactor landfill counts several unit structures - named alveoli - as shown
in the 3D rendering of a facility in Drou\`es, France (Fig. \ref{fig:droues}).
We focus on a single alveolus and we model it as a homogeneous porous medium
in which the bulk material represents the solid waste whereas the void parts
among the organic material are filled by a mixture of gases - mainly methane,
carbon dioxide, oxygen, nitrogen and water vapor - and leachates, that is a
liquid suspension based on water.
For the rest of this paper, we will refer to our domain of interest by using
indifferently the term alveolus and bioreactor, though the latter one is not
rigorous from a modeling point of view.

Each alveolus is filled with several layers of biodegradable waste and the
structure is equipped with a network of horizontal water injectors and
production pipes respectively to allow the recirculation of leachates and to
extract the gases generated by the chemical reactions. Moreover, each alveolus
is isolated from the surrounding ground in order to prevent pollutant leaks and
is covered by means of an active geomembrane. Figure \ref{fig:alveolus_scheme}
provides a schematic of an alveolus in which the horizontal pipes are organized
in order to subdivide the structure in a cartesian-like way.
Technical details on the construction and management of a bioreactor landfill
are available in \cite{itrc, ademe-fnade}.

\subsection{Physical and chemical phenomena}

Let us define the porosity $\phi$ as the fraction of void space inside
the bulk material:
\begin{equation}
\phi=\frac{\textrm{Pore Volume}}{\textrm{Total Volume}} .
\label{eq:PorosityDefinition}
\end{equation}
For biodegradable waste, we consider $\phi = 0.3$ as by experimental
measurements in \cite{Sidhu2007}, whereas it is known in the literature that
for generic waste the value drops to $0.1$.
Within this porous environment, the following phenomena take place:
\begin{itemize}
\item chemical reaction for the methane fermentation;
\item heat transfer driven by the chemical reaction;
\item transport phenomena of the gases;
\item transport and diffusion phenomena of the leachates;
\end{itemize}

Here, we briefly provide some details about the chemical reaction for the
methane generation, whereas we refer to section \ref{ref:math_model} for the
description of the remaining phenomena and the  derivation of the full
mathematical model for the coupled system.
As previously mentioned, at the beginning of the anaerobic stage the bacteria
break the previously generated molecules down to organic acids, like the
propionic acid $CH_3CH_2COOH$. This acid acts as a reacting term in the following
reaction:
\begin{equation}
CH_3CH_2COOH + H_2O \longrightarrow 3\,H_2O + CO_2 + 2CH_4.
\label{eq:chemical_reaction}
\end{equation}
The microbiota activity drives the generation of methane ($CH_4$) and is
responsible for the production of other by-products, mainly water ($H_2O$) and
carbon dioxide ($CO_2$).
As per equation (\ref{eq:chemical_reaction}), for each consumed mole of
propionic acid - equivalently referred to as organic
carbon with an abuse of notation -, two moles of methane are generated and
three moles of water and one of
carbon dioxide are produced as well.
\begin{rmrk}
We remark that in order for reaction (\ref{eq:chemical_reaction}) to take
place, water has to be added to the propionic acid. This means that the
bacteria can properly catalyze the reaction only if certain conditions on the
temperature and the humidity of the waste are fulfilled.
This paper presents a first attempt to provide a mathematical model of a
bioreactor landfill, thus both the temperature and the quantity of water inside
the facility will act as unknowns in the model (cf. section
\ref{ref:math_model}). In a more general framework, the proposed model will be
used to perform long-time forecasts of the methane generation process and the
temperature and water quantity will have the key role of control variables of
the system.
\label{rmk:humidity}
\end{rmrk}

\section{Mathematical model of a bioreactor landfill}
\label{ref:math_model}

In this section, we describe the equations modeling the phenomena taking place 
inside an alveolus. 
As stated in the introduction, the final goal of the SiViBiR++ project is 
to control and optimize the functioning of a real bioreactor landfill, hence 
a simple model to account for the phenomena under analysis is sought.
Within this framework, in this article we propose a first mathematical model 
to describe the coupled physical and chemical phenomena involved in the 
methanogenic fermentation.
In the following sections, we will provide a detailed description of the 
chemical reaction catalyzed by the methanogenic bacteria, the evolution of the 
temperature inside the alveolus and the transport phenomena driven by the 
dynamic of a mixture of gases and by the liquid water.
An extremely important aspect of the proposed model is the interaction among 
the variables at play and consequently the coupling among the corresponding 
equations.
In order to reduce the complexity of the model and to keep the corresponding 
implementation in \Feel as simple as possible, some physical phenomena have 
been neglected.
In the following subsections, we will detail the simplifying assumptions 
that allow to neglect some specific phenomena without degrading the 
reliability of the resulting model, by highlighting their limited impact 
on the global behavior of the overall system.

Let $\Omega \subset \mathbb{R}^3$ be an alveolus inside the landfill under
analysis. We split the boundary $\partial\Omega$ of the computational domain
into three non-void and non-overlapping regions $\GammaTop$, $\GammaBase$ and
$\GammaLat$, representing respectively the membrane covering the
top surface of the bioreactor, the base of the alveolus and the ground
surrounding the lateral surface of the structure.

\subsection{Consumption of the organic carbon}

As previously stated, the functioning of a bioreactor landfill relies on the
consumption of biodegradable waste by means of bacteria. From the chemical
reaction in (\ref{eq:chemical_reaction}), we may derive a relationship between
the concentration of bacteria $b$ and the concentration of the consumable
organic material which we denote by $C^{org}$. \\
The activity of the bacteria takes place if some environmental conditions are
fulfilled, namely the waste humidity and the bioreactor temperature.
Let $\wMax$ and $\Topt$ be respectively the maximal quantity of water and
the optimal temperature that allow the microbiota to catalyze the chemical
reaction (\ref{eq:chemical_reaction}). We introduce the following functions to
model the metabolism of the microbiota:
\begin{equation}
\Psi_1(w) = w \ \max \left(0, 1-\frac{w}{\wMax} \right)
,\qquad
\Psi_2(T) = \max \left(0, 1-\frac{|T-\Topt|}{A_T} \right)
\label{eq:psi1psi2}
\end{equation}
where $A_T$ is the amplitude of the variation of the temperature tolerated by the
bacteria. On the one hand, $\Psi_1$ models the fact that the bacterial
activity is proportional to the quantity of liquid water - namely leachates -
inside the bioreactor and it is prevented when the alveolus is flooded. On the
other hand, according to $\Psi_2$ the microbiota metabolism is maximum when the
current temperature equals $\Topt$ and it stops when it exceeds the interval 
of admissible temperatures $[\Topt - A_T; \Topt +A_T]$.

Since the activity of the bacteria mainly consists in consuming the organic
waste to perform reaction (\ref{eq:chemical_reaction}), it is straightforward
to deduce a proportionality relationship between $b$ and $\Corg$. By combining
the information in (\ref{eq:psi1psi2}) with this relationship, we may derive
the following law to describe the evolution of the concentration of bacteria inside
the bioreactor:
\begin{equation}
\partial_t b \propto b \ \Corg \ \Psi_1(w) \ \Psi_2(T)
\label{eq:proportionalityB}
\end{equation}
and consequently, we get a proportionality relationship for the consumption of
the biodegradable material $\Corg$:
\begin{equation}
\partial_t \Corg \propto -\partial_t b .
\label{eq:proportionalityC}
\end{equation}
Let $a_b$ and $c_b$ be two proportionality constants associated respectively with \eqref{eq:proportionalityB} and \eqref{eq:proportionalityC}.
By integrating \eqref{eq:proportionalityC} in time and introducing the proportionality constant $c_b$, we get that the concentration of bacteria reads as
\begin{equation}
b(x,t) = b_0 + c_b (\Corg_0 - \Corg(x,t))
\label{eq:bacteria}
\end{equation}
where $b_0 \coloneqq b(\cdot,0)$ and $\Corg_0 \coloneqq \Corg(\cdot,0)$ are
the initial concentrations respectively of bacteria and organic material inside the alveolus. Thus, by
plugging \eqref{eq:bacteria} into (\ref{eq:proportionalityB}) we get the
following equation for the consumption of organic carbon between the instant
$t=0$ and the final time $\Sfin$:
\begin{equation}
\begin{cases}
(1-\phi) \partial_t \Corg(x,t) = - a_b \ b(x,t) \ \Corg(x,t) \ \Psi_1(w(x,t)) \
\Psi_2(T(x,t))
\ , \ & \text{in }\Omega \times (0,\Sfin] \\
 \Corg(\cdot,0) = \Corg_0 \ , \ & \text{in }\overline{\Omega}
\end{cases}
\label{eq:organic_carbon_eq}
\end{equation}
We remark that the organic material filling the bioreactor is only present in
the bulk part of the porous medium and this is modeled by the factor $1-\phi$
which features the information about the porosity of the environment.
Moreover, we highlight that in equation (\ref{eq:organic_carbon_eq}) a
non-linear reaction term appears and in section \ref{ref:numerical_model} we
will discuss a strategy to deal with this non-linearity when moving to
the Finite Element discretization. \\
For the sake of readability, from now on we will omit the dependency on the
space and time variables in the notation for both the organic carbon and the concentration of bacteria.

\subsection{Evolution of the temperature}

The equation describing the evolution of the temperature $T$ inside the
bioreactor is the classical heat equation with a source term proportional to
the consumption of bacteria.
We consider the external temperature to be fixed by imposing Dirichlet boundary
conditions on $\partial\Omega$.
\begin{rmrk}
Since we are interested in the long-time evolution of the system
($\Sfin =40$ years), the unit time interval is
sufficiently large to allow daily variations of the temperature to be neglected.
Moreover, we assume that the external temperature remains constant during the
whole life of the bioreactor. From a physical point of view, this assumption is
not realistic but we conjecture that only small fluctuations would arise by the
relaxation of this hypothesis. A future improvement of the model may focus on
the integration of dynamic boundary conditions in order to model seasonal
changes of the external temperature.
\end{rmrk}
The resulting equation for the temperature reads as follows:
\begin{equation}
\begin{cases}
 \partial_t T(x,t) - k_T\Delta T(x,t) = -c_T \partial_t \Corg(x,t)  \ ,
\ & \text{in} \ \Omega \times (0,\Sfin] \\
 T(x,\cdot) = \Tm \ , \ & \text{on} \ \GammaTop \times
(0,\Sfin] \\
 T(x,\cdot) = \Tg \ , \ & \text{on} \ \GammaBase \cup \GammaLat
\times (0,\Sfin]
\\
 T(\cdot,0) = T_0 \ , \ & \text{in} \ \overline{\Omega}
\end{cases}
\label{eq:heat_eq}
\end{equation}
where $k_T$ is the thermal conductivity of the biodegradable waste and
$c_T$ is a scaling factor that accounts for the heat transfer due to the
chemical reaction catalyzed by the bacteria. The values $\Tm, \ \Tg \
\text{and} \ T_0$ respectively represent the external temperature on the
membrane $\GammaTop$, the external temperature of the ground
$\GammaBase \cup \GammaLat$ and the initial temperature inside the
bioreactor. \\
For the sake of readability, from now on we will omit the dependency on the
space and time variables in the notation of the temperature.

\subsection{Velocity field of the gas}
\label{ref:Darcy}

In order to model the velocity field of the gas inside the bioreactor, we have
to introduce some assumptions on the physics of the problem. First of all, we
assume the gas to be incompressible. This hypothesis stands if a very slow
evolution of the mixture of gases takes place and this is the case for a
bioreactor landfill in which the methane fermentation gradually decreases along
the 40 years lifetime of the facility. 
Additionally, the decompression generated by the extraction of the gases through 
the pipes is negligible due to the weak gradient of pressure applied to the 
production system.
Furthermore, we assume low Reynolds and low Mach numbers for the problem under 
analysis: this reduces to having a laminar slow flow which, as previously stated, 
is indeed the dynamic taking place inside an alveolus. 
Eventually, we neglect the effect due to the gravity on the dynamic of the 
mixture of gases: owing to the small height of the alveolus (approximately 
$\SI{90}{m}$), the variation of the pressure in the vertical direction due to 
the gravity is 
limited and in our model we simplify the evolution of the gas by neglecting 
the hydrostatic component of the pressure.

Under the previous assumptions, the behavior of the gas mixture inside 
a bioreactor landfill may be 
described by a mass balance equation coupled with a Darcy's law
\begin{equation}
\begin{cases}
\nabla \cdot u = 0 \ , \ & \text{in} \ \Omega \\
u = - \nabla p \ , \ & \text{in} \ \Omega
\end{cases}
\label{eq:darcy}
\end{equation}
where $p \coloneqq \frac{D}{\phi \mu_\text{gas}} P$, $D$ is the permeability of
the porous medium, $\phi$ its porosity and $\mu_\text{gas}$ the gas viscosity
whereas $P$ is the pressure inside the bioreactor. 
In (\ref{eq:darcy}), the incompressibility assumption has been expressed by 
stating that the gas flow is isochoric, that is the velocity is divergence-free.
This equation is widely used in the literature to model porous media (cf. e.g. 
\cite{cordoba2007analytical, lebeau2009natural}) and provides a coherent 
description of the phenomenon under analysis in the bioreactor landfill. 
As a matter of fact, it is reasonable to assume that the density of the gas 
mixture is nearly constant inside the domain, owing to the weak gradient of pressure 
applied to extract the gas via the production system and to the slow rate of 
methane generation via the fermentation process, that lasts approximately 40 
years.

To fully describe the velocity field, the effect of the production system that
extracts the gases from the bioreactor has to be accounted for.
We model the production system as a set of $\Nprod$ cylinders $\Prod^i$'s thus
the effect of the gas extraction on each pipe results in a condition on the
outgoing flow.
Let $\Jout > 0$ be the mass flow rate exiting from the alveolus through each
production pipe. The system of equations (\ref{eq:darcy}) is coupled with the
following conditions on the outgoing normal flow on each drain used to
extract the gas:
\begin{equation}
\int_{(\partial\Prod^i)^{\textrm{n}}} (\Cdx+M+O+N+h)\, u \cdot n \, d\sigma =
\Jout \; \qquad \forall i=1,\ldots,\Nprod.
\label{eq:productors}
\end{equation}
In (\ref{eq:productors}), $n$ is the outward
normal vector to the surface, $(\partial\Prod^i)^{\textrm{n}}$ is the part of
the boundary of the cylinder $\Prod^i$ which belongs to the lateral surface of
the alveolus itself and the term $(\Cdx+M+O+N+h)$ represents the total concentration
of the gas mixture starring carbon dioxide, methane, oxygen, nitrogen and water
vapor. \\
Since the cross sectional area of the pipes belonging to the production system
is
negligible with respect to the size of the overall alveolus, we model these
drains as 1D lines embedded in the 3D domain. Owing to this, in the
following subsection we present a procedure to integrate the information
(\ref{eq:productors}) into a source term named $\Fout$ in order to simplify
the problem that describes the dynamic of the velocity field inside a
bioreactor landfill.
\begin{rmrk}
According to conditions (\ref{eq:productors}), the velocity $u$ depends on the
concentrations of the gases inside the bioreactor, thus is a function of both
space and time. Nevertheless, the velocity field at each time step is
independent from the previous ones and is only influenced by the distribution
of gases inside the alveolus. For this reason, we neglect the dependency on the
time variable and we consider $u$ being only a function of space.
\end{rmrk}

\subsubsection{The source term $\Fout$}
\label{ref:Fout}

As previously stated, each pipe $\Prod^i$ is modeled as a cylinder of radius
$R$ and length $L$. Hence, the cross sectional area
$(\partial\Prod^i)^{\textrm{n}}$ and the lateral surface
$(\partial\Prod^i)^{\textrm{l}}$ respectively measure $\pi R^2$ and $2\pi R L$.
\\
We assume the gas inside the cylinder to instantaneously exit the alveolus
through its boundary $(\partial\Prod^i)^{\textrm{n}}$, that is the outgoing
flow (\ref{eq:productors}) is equal to the flow entering the drain through its
lateral surface. Thus we may neglect the gas dynamic inside the pipe and
(\ref{eq:productors}) may be rewritten as
\begin{equation}
\int_{(\partial\Prod^i)^{\textrm{n}}} (\Cdx+M+O+N+h)\, u \cdot n \, d\sigma =
\int_{(\partial\Prod^i)^{\textrm{l}}} (\Cdx+M+O+N+h)\, u \cdot n \, d\sigma =
\Jout \qquad \forall i=1,\ldots,\Nprod.
\label{eq:productorsLAT}
\end{equation}
Moreover, under the hypothesis that the quantity of gas flowing from
the bioreactor to the inside of the cylinder $\Prod^i$ is uniform over its
lateral surface, that is the same gas mixture surrounds the drain in all the
points along its dominant size, we get
\begin{equation}
\left( \int_{(\partial\Prod^i)^{\textrm{l}}} (\Cdx+M+O+N+h)\, d\sigma \right) u
\cdot n = \Jout \,\; \text{on} \,\; (\partial\Prod^i)^{\textrm{l}} \qquad
\forall i=1,\ldots,\Nprod.
\label{eq:lateralBC}
\end{equation}
We remark that gas densities may be considered uniform along the perimeter of
the cylinder only if the latter is small enough, that is the aforementioned
assumption is likely to be true if the radius of the pipe is small in
comparison with the size of the alveolus. Within this framework,
(\ref{eq:lateralBC}) reduces to
\begin{equation}
u \cdot n = \frac{\Jout}{\displaystyle 2\pi R \int_{\Line^i} (\Cdx+M+O+N+h)\,
dl} \,\; \text{on} \,\; (\partial\Prod^i)^{\textrm{l}} \qquad
\forall i=1,\ldots,\Nprod
\label{eq:limitBC}
\end{equation}
where $\Line^i$ is the centerline associated with the cylinder $\Prod^i$.
By coupling (\ref{eq:darcy}) with (\ref{eq:limitBC}), we get the following PDE
to model the velocity field:
\begin{equation}
\begin{cases}
- \Delta p = 0 \ , \ & \text{in} \ \Omega \\
\nabla p \cdot n = -\frac{\displaystyle \Jout}{\displaystyle 2\pi R
\int_{\Line^i} (\Cdx+M+O+N+h)\,
dl} \ , \ & \text{on} \ (\partial\Prod^i)^{\textrm{l}} \qquad
\forall i=1,\ldots,\Nprod
\end{cases}
\label{eq:velocityField}
\end{equation}
Let us consider the variational formulation of problem
(\ref{eq:velocityField}): we seek $p \in H^1(\Omega)$ such that
\begin{equation}
\int_\Omega \nabla p \cdot \nabla \delta p \, dx = \sum_{i=1}^{\Nprod}
-\frac{\displaystyle \Jout}{\displaystyle 2\pi R \int_{\Line^i}
(\Cdx+M+O+N+h)\,
dl} \, \int_{(\partial\Prod^i)^{\textrm{l}}} \delta p \, d\sigma \qquad \forall
\delta p \in \mathcal{C}^1_0(\Omega).
\label{eq:variationalVel}
\end{equation}
We may introduce the term $\Fout$ as the limit when $R$ tends to zero of
the right-hand side of (\ref{eq:variationalVel}):
\begin{equation}
\Fout \coloneqq \sum_{i=1}^{\Nprod} -\frac{\displaystyle \Jout}{\displaystyle
\int_{\Line^i} (\Cdx+M+O+N+h)\, dl} \delta_{\Line^i}
\label{eq:Fout}
\end{equation}
where $\delta_{\Line^i}$ is a Dirac mass concentrated along the centerline
$\Line^i$ of the pipe $\Prod^i$. \\
Hence, the system of equations describing the evolution of the velocity inside
the alveolus may be written as
\begin{equation}
\begin{cases}
\nabla \cdot u = \Fout \ , \ & \text{in} \ \Omega \\
u = - \nabla p \ , \ & \text{in} \ \Omega \\
u \cdot n = 0 \ , \ & \text{on} \ \partial\Omega
\end{cases}
\label{eq:velocity1D}
\end{equation}
where the right-hand side of the mass balance equation may be
either (\ref{eq:Fout}) or a mollification of it.

\subsection{Transport phenomena for the gas components}
\label{ref:gases}

Inside a bioreactor landfill the pressure field is comparable to the external 
atmospheric pressure.
This low-pressure does not provide the physical conditions for gases to liquefy. 
Hence, the gases are not present in liquid phase and solely the dynamic of the gas 
phases has to be accounted for.
Within this framework, 
in section \ref{ref:vaporLiquid} we consider the case of water for which phase 
transitions driven by heat transfer phenomena are possible, 
whereas in the current section we focus on the remaining gases (i.e. oxygen, nitrogen, 
methane and carbon dioxide) which solely exist in gas phase.

Let $u$ be the velocity of the gas mixture inside the alveolus. We
consider a generic gas whose concentration inside the bioreactor is named $G$.
The evolution of $G$ fulfills the classical pure advection equation:
\begin{equation}
\begin{cases}
\phi \partial_t G(x,t) + u \cdot \nabla G(x,t) = F^G(x,t)  \ ,
\ & \text{in} \ \Omega \times (0,\Sfin] \\
G(\cdot,0) = G_0 \ , \ & \text{in} \ \overline{\Omega}
\end{cases}
\label{eq:advectionG}
\end{equation}
where $\phi$ is again the porosity of the waste.
The source term $F^G(x,t)$ depends on the gas and will be detailed in the
following subsections.

\subsubsection{The case of oxygen and nitrogen}

We recall that the oxygen concentration is named $O$, whereas the nitrogen one
is $N$.
Neither of these components appears in reaction (\ref{eq:chemical_reaction})
thus the associated source terms are $F^O(x,t)=F^N(x,t)=0$. The resulting
equations (\ref{eq:OxyNitro}) are closed by the initial conditions $O(\cdot,0)
=
O_0$ and $N(\cdot,0) = N_0$.
\begin{equation}
  \begin{split}
    &\phi \partial_t O + u \cdot \nabla O = 0 \\
    &\phi \partial_t N + u \cdot \nabla N = 0 \\
  \end{split}
\label{eq:OxyNitro}
\end{equation}
Both the oxygen and the nitrogen are extracted by the production system thus
their overall concentration may be negligible with respect to the quantity of
carbon dioxide and methane inside the alveolus. Hence, for the rest of this
paper we will neglect equations (\ref{eq:OxyNitro}) by considering $O(x,t)
\simeq O_0 \simeq 0$ and $N(x,t) \simeq N_0 \simeq 0$.

\subsubsection{The case of methane and carbon dioxide}

As previously stated, (\ref{eq:chemical_reaction}) describes the methanogenic
fermentation that starting from the propionic acid drives the production of
methane, having carbon dioxide as by-product.
Equation (\ref{eq:advectionG}) stands for both the methane $M$ and the carbon
dioxide $\Cdx$. For these components, the source terms have to account for the
production of gas starting from the transformation of biodegradable waste.
Thus, the source terms are proportional to the consumption of the quantity
$\Corg$ through some constants $c_M$ and $c_C$ specific to the chemical
reaction and the component:
$$
F^j(x,t) = -c_j \partial_t \Corg \quad , \quad j=M,C
$$
In a similar fashion as before, the resulting equations read as
\begin{equation}
  \begin{split}
    &\phi \partial_t M + u \cdot \nabla M = -c_M \partial_t \Corg \\
    &\phi \partial_t \Cdx + u \cdot \nabla \Cdx = -c_C \partial_t \Corg \\
  \end{split}
\label{eq:MethaneCdx}
\end{equation}
and they are coupled with appropriate initial conditions $M(\cdot,0) = M_0$ and
$\Cdx(\cdot,0) = \Cdx_0$.

\subsection{Dynamic of water vapor and liquid water}
\label{ref:vaporLiquid}

Inside a bioreactor landfill, water exists both in vapor and liquid phase. 
Let $h$ be the concentration of water vapor and $w$ the one of liquid water.
The variation of temperature responsible for phase transitions inside the 
alveolus is limited, whence we do not consider a two-phase flow for the water 
but we describe separately the dynamics of the gas and liquid phases of the fluid. 
On the one hand, the water vapor inside the bioreactor landfill evolves as the 
gases presented in section \ref{ref:gases}: it is produced by the chemical reaction 
(\ref{eq:chemical_reaction}), it is transported by the velocity 
field $u$ and is extracted via the pipes of the production system;
as previously stated, no effect of the gravity is accounted for. 
This results in a pure advection equation for $h$. 
On the other hand, the dynamic of the liquid water may be schematized as follows: 
it flows in through the injector system at different levels of the alveolus, 
is transported by a vertical field $u_w$ due to the effect of gravity and 
is spread within the porous medium. 
The resulting governing equation for $w$ is an advection-diffusion equation.
Eventually, the phases $h$ and $w$ are coupled by a source term that accounts for 
phase transitions.

Owing to the different nature of the phenomena under analysis and to the limited 
rate of heat transfer inside a bioreactor landfill, in the rest of this section we 
will describe separately the equations associated with the dynamics of the water vapor 
and the liquid water, highlighting their coupling due to the phase transition 
phenomena.

\subsubsection{Phase transitions}
\label{ref:phase}

Two main phenomena are responsible for the production of water vapor inside a
bioreactor landfill. On the one hand, vapor is a product of the chemical
reaction (\ref{eq:chemical_reaction}) catalyzed by the microbiota during the
methanogenic fermentation process. On the other hand, heat transfer causes part
of the water vapor to condensate and part of the liquid water to evaporate. \\
Let us define the vapor pressure of water $\Pvap$ inside the alveolus as the
pressure at which water vapor is in thermodynamic equilibrium with its
condensed state. Above this critical pressure, water vapor condenses, that is
it turns to the liquid phase. This pressure is proportional to the
temperature $T$ and may be approximated by the following Rankine law:
\begin{equation}
\Pvap(T) = P_0\exp \left( s_0 - \frac{s_1}{T} \right)
\label{eq:vapor_pressure}
\end{equation}
where $P_0$ is a reference pressure, $s_0$ and $s_1$ are two constants known
by experimental results and $T$ is the temperature measured in Kelvin.
If we restrict to a range of moderate temperatures, we can approximate the
exponential in (\ref{eq:vapor_pressure}) by means of a linear law. Let
$H_0$ and $H_1$ be two known constants, we get
\begin{equation}
\Pvap(T) \simeq H_0 + H_1 T .
\label{eq:linearExp}
\end{equation}
Let $P^h$ be the partial pressure of the water vapor inside the gas mixture. We
can compute $P^h$ multiplying the total pressure $p$ by a
scaling factor representing the ratio of water vapor inside the gas mixture:
\begin{equation}
P^h = \frac{h}{\Cdx + M + O + N + h} p .
\label{eq:partial_pressure}
\end{equation}
The phase transition process features two different phenomena. On the one hand,
when the pressure $P^h$ is higher than the vapor pressure of water $\Pvap$ the
vapor condensates.
By exploiting (\ref{eq:partial_pressure}) and (\ref{eq:linearExp}), the
condition $P^h > \Pvap(T)$ may be rewritten as
$$
h - H(T) > 0 \qquad , \qquad H(T) \coloneqq (\Cdx + M + O + N + h)\frac{H_0 +
H_1 T}{p} .
$$
We assume the phase transition to be instantaneous, thus the
condensation of water vapor may be expressed through the following function
\begin{equation}
\Fcond \coloneqq c_{h\to w}\max(h - H(T), 0)
\label{eq:condensation}
\end{equation}
where $c_{h\to w}$ is a scaling factor. As per (\ref{eq:condensation}), the
production of vapor from liquid water is $0$ as soon as the concentration of vapor is
larger than the threshold $H(T)$, that is the air is saturated. \\
In a similar fashion, we may model the evaporation of liquid water. When $P^h$
is below the vapor pressure of water $\Pvap$ - that is $h - H(T) < 0$ - part of
the liquid water generates vapor. The evaporation rate is proportional to the
difference $\Pvap(T) - P^h$ and to the quantity of liquid water $w$ available
inside the alveolus. Hence, the evaporation of liquid water is modeled by the
following expression
\begin{equation}
\Fevap \coloneqq c_{w\to h} \max(H(T) - h, 0) w .
\label{eq:evaporation}
\end{equation}
\begin{rmrk}
Since the quantity of water vapor inside a bioreactor landfill is
negligible, we assume that the evaporation process does not
significantly affect the dynamic of the overall system. Hence, in the rest of
this paper, we will neglect this phenomenon by modeling only the
condensation (\ref{eq:condensation}).
\end{rmrk}

\subsubsection{The case of water vapor}

The dynamic of the water vapor may be modeled using (\ref{eq:advectionG}) as
for the other gases. In this case, the source term has to account for both
the production of water vapor due to the chemical reaction
(\ref{eq:chemical_reaction}) and its decrease as a consequence of the
condensation phenomenon:
$$
F^h(x,t) \coloneqq -c_h \partial_t \Corg - \Fcond
$$
where $c_h$ is a scaling factor describing the relationship between the
consumption of organic carbon and the generation of water vapor.
The resulting advection equation reads as
\begin{equation}
\phi \partial_t h + u \cdot \nabla h = -c_h \partial_t \Corg - \Fcond
\label{eq:water_vapor}
\end{equation}
and it is coupled with the initial condition $h(\cdot,0) = h_0$.

\subsubsection{The case of liquid water}

The liquid water inside the bioreactor is modeled by an advection-diffusion
equation in which the drift term is due to the gravity, that is the transport
phenomenon is mainly directed in the vertical direction and is associated with
the liquid flowing downward inside the alveolus. 
\begin{equation}
\begin{cases}
\phi \partial_t w(x,t) + u_w \cdot \nabla w(x,t) - k_w \Delta w(x,t) =
\Fcond(x,t) \ , \ & \text{in} \ \Omega \times (0,\Sfin] \\
k_w \nabla w(x,t) \cdot n = 0 \ \text{and} \ u_w \cdot n = 0 \ , \ & \text{on}
\
\GammaTop \cup \GammaLat \times (0,\Sfin] \\
k_w \nabla w(x,t) \cdot n = 0  \ , \ & \text{on} \ \GammaBase \times (0,\Sfin]
\\
w(\cdot,0) = w_0 \ , \ & \text{in} \ \overline{\Omega}
\end{cases}
\label{eq:Liquid}
\end{equation}
where $u_w \coloneqq (0,0,-\| u_w \|)^T$ is the vertical velocity of the water
and $k_w$ is its diffusion coefficient. The right-hand side of the first
equation accounts for the water production by condensation as described in
section \ref{ref:phase}. On the one hand, the free-slip boundary conditions on
the lateral and top surfaces allow water to slide but prevent its exit, that is
the top and lateral membranes are waterproof. On the other hand, the
homogeneous Neumann boundary condition on the bottom of the domain describes
the ability of the water to flow through this membrane. These conditions are
consistent with the impermeability of the geomembranes and with the
recirculation of leachates which are extracted when they accumulate in the
bottom part of the alveolus and are reinjected in the upper
layers of the waste management facility.

\begin{rmrk}
It is well-known in the literature that the evolution of an incompressible 
fluid inside a given domain is described by the Navier-Stokes equation.
In (\ref{eq:Liquid}), we consider a simplified version of the aforementioned 
equation by linearizing the inertial term. 
As previously stated, the dynamic of the fluids inside the bioreactor landfill 
is extremely slow and we may assume a low Reynolds number regime for the water 
as well.
Under this assumption, the transfer of kinetic energy in the turbulent cascade 
due to the non-linear term of the Navier-Stokes equation may be neglected.
Moreover, by means of a linearization of the inertial term $(w \cdot \nabla ) w$, the transport effect 
is preserved and the resulting parabolic advection-diffusion problem 
(\ref{eq:Liquid}) may be interpreted as an unsteady version of the classical 
Oseen equation \cite{MR1218879}.
\end{rmrk}

\begin{rmrk}
Equation \eqref{eq:Liquid} may be furtherly interpreted as a special advection 
equation modeling the transport phenomenon within a porous medium.
As a matter of fact, the diffusion term $-k_w \Delta w$ accounts for the 
inhomogeneity of the environment in which the water flows and describes the 
fact that the liquid spreads in different directions while flowing downwards
due to the encounter of blocking solid material along its path.
The distribution of the liquid into different directions is random and is mainly 
related to the nature of the surrounding environment thus we consider an 
isotropic diffusion tensor $k_w$.
The aforementioned equation is widely used (cf. e.g. \cite{vafai2015handbook}) 
to model flows in porous media and is strongly connected with the description of the 
porous environment via the Darcy's law introduced in section \ref{ref:Darcy}. \\
Within the framework of our problem, the diffusion term is extremely important 
since it models the spread of water and leachates inside the bioreactor landfill 
and the consequent humidification of the whole alveolus and not solely of the 
areas neighboring the injection pipes. 
\end{rmrk}

Eventually, problem (\ref{eq:Liquid}) is closed by a set of
conditions that describe the injection of liquid water and leachates through
$\Ninj$ pipes $\Inj^i$'s. As previously done for the production system, we
model each injector as a cylinder of radius $R$ and length $L$ and we denote by
$(\partial\Inj^i)^{\textrm{n}}$ and $(\partial\Inj^i)^{\textrm{l}}$
respectively the part of the boundary of the cylinder which belongs to the
boundary of the bioreactor and its lateral surface. The aforementioned inlet
condition reads as
$$
\int_{(\partial\Inj^i)^{\textrm{n}}} k_w \nabla w \cdot n \, d\sigma =
-\Jin \; \qquad \forall i=1,\ldots,\Ninj
$$
where $\Jin > 0$ is the mass flow rate entering the alveolus through each
injector.
As for the production system in section \ref{ref:Fout}, we may now
integrate this condition into a source term for equation (\ref{eq:Liquid}).
Under the assumption that the flow is instantaneously distributed along the
whole cylinder in a uniform way, we get
$$
\int_{(\partial\Inj^i)^{\textrm{n}}} k_w \nabla w \cdot n \, d\sigma =
-\int_{(\partial\Inj^i)^{\textrm{l}}} k_w \nabla w \cdot n \, d\sigma =
-\Jin \qquad \forall i=1,\ldots,\Ninj.
$$
Consequently the condition on each injector reads as
$$
k_w \nabla w \cdot n = \frac{\Jin}{2\pi R L} \,\; \text{on} \,\;
(\partial\Inj^i)^{\textrm{l}} \qquad
\forall i=1,\ldots,\Ninj
$$
and we obtain the following source term $\Fin$
\begin{equation}
\Fin \coloneqq \sum_{i=1}^{\Ninj} \frac{\displaystyle \Jin}{\displaystyle L}
\delta_{\Line^i}
\label{eq:Fin}
\end{equation}
where $\delta_{\Line^i}$ is a Dirac mass concentrated along the centerline
$\Line^i$ of the pipe $\Inj^i$.
Hence, the resulting dynamic of the liquid water inside an alveolus is modeled
by the following PDE:
\begin{equation}
\begin{cases}
\phi \partial_t w + u_w \cdot \nabla w - k_w \Delta w =
F^w \ , \ & \text{in} \ \Omega \times (0,\Sfin] \\
k_w \nabla w \cdot n = 0 \ \text{and} \ u_w \cdot n = 0 \ , \ & \text{on} \
\GammaTop \cup \GammaLat \times (0,\Sfin] \\
k_w \nabla w \cdot n = 0  \ , \ & \text{on} \ \GammaBase \times (0,\Sfin] \\
\end{cases}
\label{eq:water1D}
\end{equation}
with $F^w \coloneqq \Fcond + \Fin$ and the initial condition
$w(\cdot,0) = w_0$. By analyzing the right-hand side of equation
(\ref{eq:water1D}), we remark that neglecting the effect of evaporation in the
phase transition allows to decouple the dynamics of liquid water and water
vapor. Moreover, as previously stated for equation (\ref{eq:velocity1D}),
$\Fin$
may be chosen either according to definition (\ref{eq:Fin}) or by means of an
appropriate mollification.

\section{Numerical approximation of the coupled system}
\label{ref:numerical_model}

This section is devoted to the description of the numerical strategies used to
discretize the fully coupled model of the bioreactor landfill introduced in
section \ref{ref:math_model}.
We highlight that one of the main difficulties of the presented model is the
coupling of all the equations and the multiphysics nature of the problem under
analysis.
Here we propose a first attempt to discretize the full model by introducing an
explicit coupling of the equations, that is by considering the source term in
each equation as function of the variables at the previous iteration.

\subsection{Geometrical model of an alveolus}

\begin{figure}[htbp]
\centering
\subfloat[Cubic reference domain.]
{
\includegraphics[width=0.31\columnwidth]{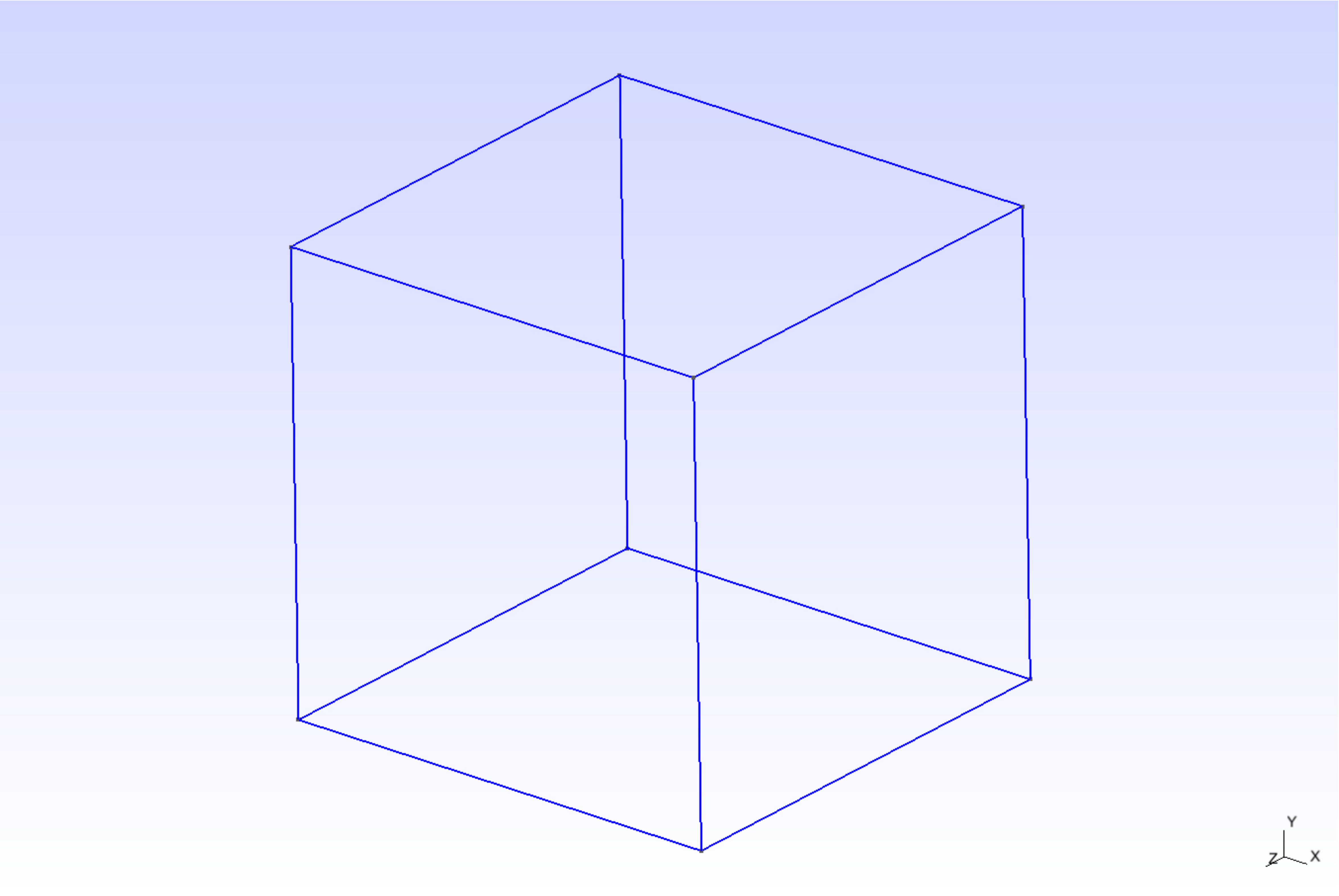}
\label{fig:cube}
}
\hfill
\subfloat[Lateral shear on a face.]
{
\includegraphics[width=0.31\columnwidth]{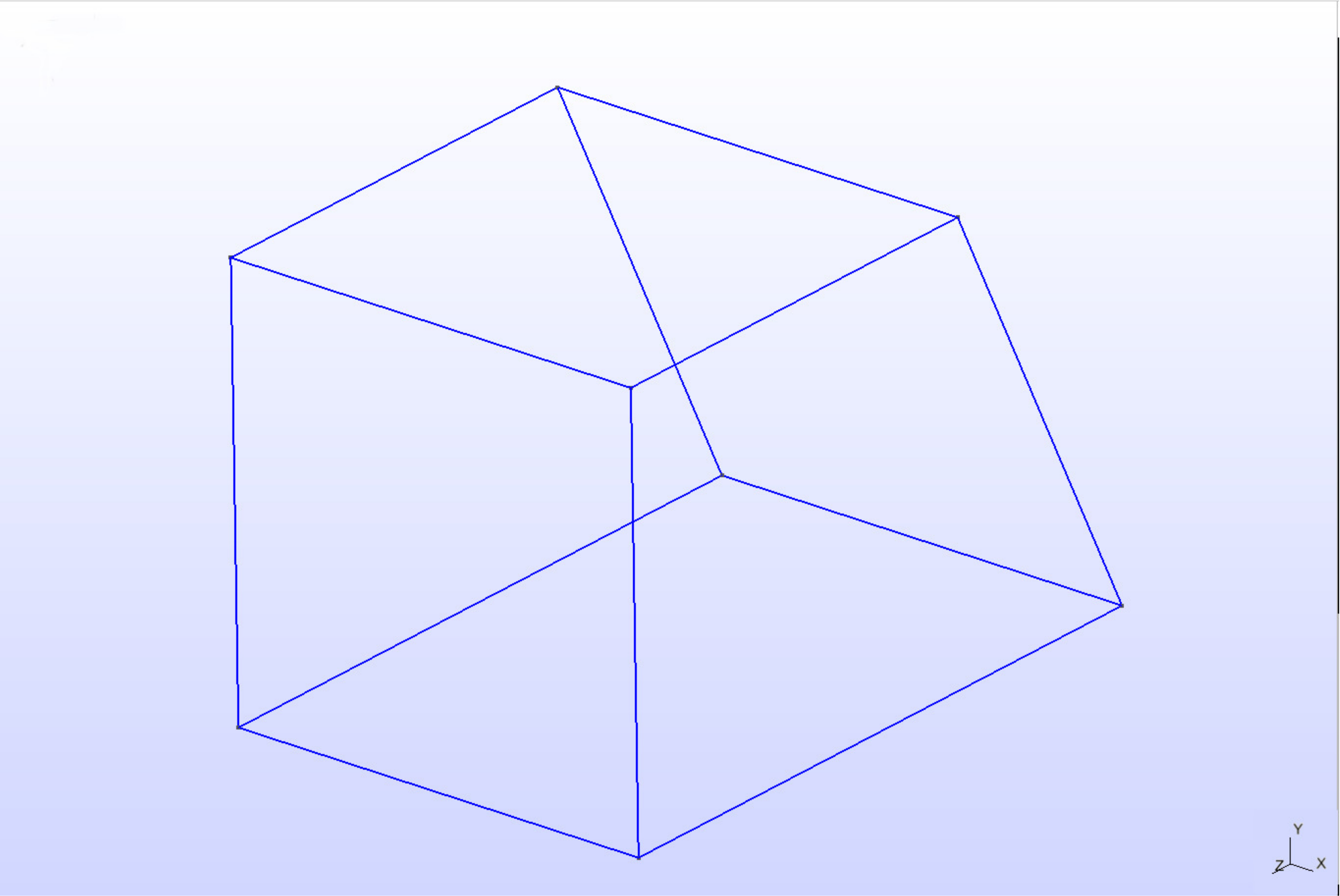}
\label{fig:lateral_shear}
}
\hfill
\subfloat[Pyramidal domain.]
{
\includegraphics[width=0.31\columnwidth]{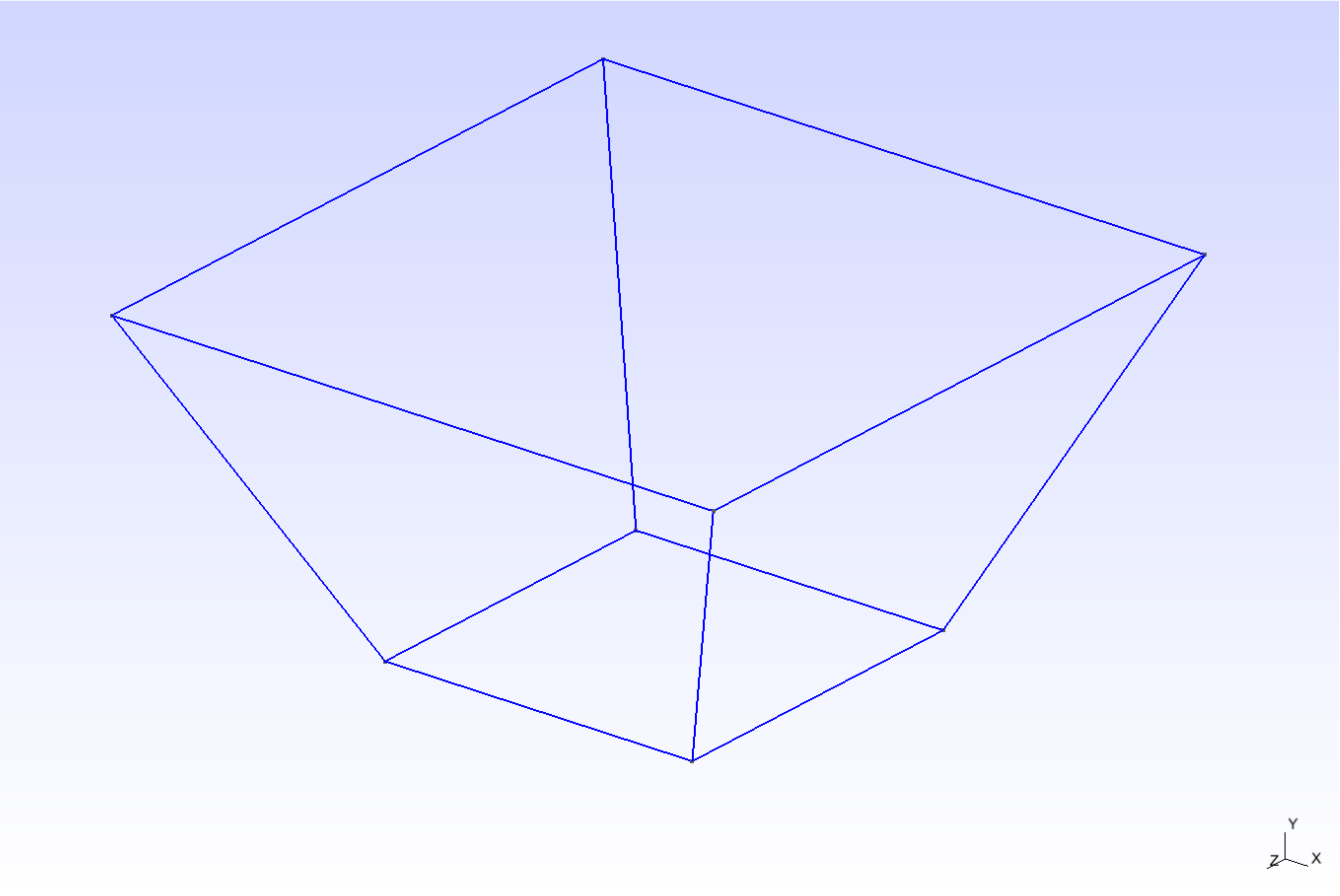}
\label{fig:pyramid}
}
\caption{Reference domain for an alveolus and admissible transformations.}
\label{fig:gmsh_model}
\end{figure}

As previously stated, a bioreactor landfill is composed by several alveoli.
Each alveolus may be modeled as an independent structure obtained starting from
a cubic reference domain (Fig. \ref{fig:cube}) to which pure shear
transformations are applied (Fig. \ref{fig:lateral_shear}-\ref{fig:pyramid}).
For example, the pure lateral shear in figure \ref{fig:lateral_shear} allows to
model the left-hand side alveolus in figure \ref{fig:alveolus_scheme} whereas
the right-hand side one may be geometrically approximated by means of the
pyramid in figure \ref{fig:pyramid}. \\
In order to model the network of the water injectors and the one of the drains
extracting the gas, the geometrical domains in figure \ref{fig:gmsh_model} are
equipped with a cartesian distribution of horizontal lines, the 1D model being
justified by the assumption in section \ref{ref:math_model}.

\subsection{Finite Element approximations of the organic carbon and heat equations}
\label{ssec:num_organic_heat_eq}

Both equations (\ref{eq:organic_carbon_eq}) and (\ref{eq:heat_eq}) are
discretized using Lagrangian Finite Element functions. In particular, the time
derivative is approximated by means of an implicit Euler scheme, whereas the
basis functions for the spatial discretization are the classical $\mathbb{P}_k$
Finite Element functions of degree $k$. \\
Let $t=t_n$. We consider the following quantities at time $t_n$ as known
variables: $\Corg_n \coloneqq \Corg(x,t_n)$, $T_n \coloneqq T(x,t_n)$ and $w_n
\coloneqq w(x,t_n)$.
The consumption of organic carbon is described by equation
(\ref{eq:organic_carbon_eq}) coupled with equation (\ref{eq:bacteria}) for the
dynamic of the bacteria. At each time step, we seek $\Corg_{n+1} \in
H^1(\Omega)$ such that
$$
\int_\Omega (1-\phi) \frac{\Corg_{n+1} - \Corg_n}{\Delta t} \delta C \, dx = -
\int_\Omega a_b \ \Corg_{n+1} \left[ b_0 + c_b
(\Corg_0 - \Corg_n) \right] \Psi_1(w_n) \
\Psi_2(T_n) \ \delta C \, dx \qquad \forall \delta C \in H^1(\Omega)
$$
We remark that in the previous equation the non-linear reaction term has been
handled in a semi-implicit way by substituting $(\Corg_{n+1})^2$ by
$\Corg_{n+1}\Corg_n$ in the right-hand side.
Hence, the bilinear and linear forms associated with the variational
formulation at $t=t_n$ respectively read as
\begin{equation}
a_{\Corg}(\Corg_{n+1},\delta C) = \int_\Omega A_C \ \Corg_{n+1} \ \delta C \, dx
\qquad , \qquad
l_{\Corg}(\delta C) = \int_\Omega (1-\phi) \ \Corg_n \ \delta C \, dx
\label{eq:variational_Corg}
\end{equation}
where $A_C = (1-\phi) + \Delta t \ a_b \left[ b_0 + c_b
(\Corg_0 - \Corg_n) \right] \Psi_1(w_n) \ \Psi_2(T_n)$ and
$a_{\Corg}(\Corg_{n+1},\delta C) = l_{\Corg}(\delta C) \ \forall \delta C \in
H^1(\Omega)$.

In a similar fashion, we derive the variational formulation of the heat
equation and at $t=t_n$ we seek $T_{n+1} \in H^1(\Omega)$ such that
$T_{n+1}|_{\GammaTop} = T_{\textrm{m}}$, $T_{n+1}|_{\GammaBase \cup
\GammaLat} = T_{\textrm{g}}$ and $a_T(T_{n+1},\delta T) = l_T(\delta
T) \ \forall \delta T \in H^1_0(\Omega)$, where
\begin{gather}
a_T(T_{n+1},\delta T) = \int_\Omega \ T_{n+1} \ \delta T \, dx +
\Delta t \int_\Omega \ k_T \nabla T_{n+1} \cdot \nabla \delta T \, dx \, ,
\label{eq:variational_a_T} \\
l_T(\delta T) = \int_\Omega \ T_n \ \delta T \, dx - \int_\Omega c_T
(\Corg_{n+1}-\Corg_n) \ \delta T \, dx .
\label{eq:variational_l_T}
\end{gather}
\begin{rmrk}
In (\ref{eq:variational_a_T}) and (\ref{eq:variational_l_T}), we evaluate the
time derivative $\partial_t \Corg$ at $t=t_n$, that is we consider the
current value and not the previous one as stated at the beginning of this
section. This is feasible because the solution of problem
(\ref{eq:variational_Corg}) precedes the one of the heat equation thus the value
$\Corg_{n+1}$ is known when solving
(\ref{eq:variational_a_T})-(\ref{eq:variational_l_T}).
\end{rmrk}
\begin{rmrk}
In (\ref{eq:variational_l_T}) we
assume that the same time discretization is used for both the organic carbon
and the temperature. If this is not the case, the second term in the
linear form $l_T(\cdot)$ would feature a scaling factor $\frac{\Delta
t_T}{\Delta t_C}$, the numerator being the time scale associated with the
temperature and the denominator the one for the organic carbon. For the rest of
this paper, we will assume that all the unknowns are approximated using the
same time discretization.
\label{rmk:deltaT}
\end{rmrk}
By substituting $\Corg_{n}$ and $T_n$ with their Finite Element
counterparts $\Corg_{h,n}$ and $T_{h,n}$ in (\ref{eq:variational_Corg}),
(\ref{eq:variational_a_T}) and (\ref{eq:variational_l_T}) we obtain the
corresponding discretized equations for the organic carbon and the temperature.

\subsection{Stabilized dual-mixed formulation of the velocity field}

A good approximation of the velocity field is a key point for a satisfactory
simulation of all the transport phenomena. In order for problem
(\ref{eq:velocity1D}) to be well-posed, the following
compatibility
condition has to be fulfilled
\begin{equation}
\int_\Omega \Fout = 0.
\label{eq:compatibility}
\end{equation}
Nevertheless, (\ref{eq:compatibility}) does not stand for the problem under
analysis thus we consider a slightly modified version of problem
(\ref{eq:velocity1D}) by introducing a small perturbation parameter $\lambda =
\mathcal{O}(\ell_K)$, $\ell_K$ being the diameter of the element $K$ of the
triangulation $\mathcal{T}_h$:
\begin{equation}
\begin{cases}
\nabla \cdot u + \lambda p = \Fout \ , \ & \text{in} \ \Omega \\
u = - \nabla p \ , \ & \text{in} \ \Omega \\
u \cdot n = 0 \ , \ & \text{on} \ \partial\Omega \\
\end{cases}
\label{eq:vel_perturbed}
\end{equation}
Hence the resulting problem (\ref{eq:vel_perturbed}) is well-posed even if
(\ref{eq:compatibility}) does not stand.

It is well-known in the literature \cite{WRCR:WRCR6472} that classical
discretizations of problem (\ref{eq:vel_perturbed}) by means of Lagrangian
Finite Element functions lead to poor approximations of the velocity field. A
widely accepted workaround relies on the derivation of mixed formulations in
which a simultaneous approximation of pressure and velocity fields is
performed by using different Finite Element spaces \cite{RaviartThomas}.

\subsubsection{Dual-mixed formulation}

Let $\HDiv = \{ v \in [L^2(\Omega)]^3 \ \text{s.t.} \ \nabla \cdot v \in
L^2(\Omega) \}$ and $\Hdiv0 = \{ v \in \HDiv \ \text{s.t.} \ v \cdot n = 0 \
\text{on} \ \partial\Omega \}$. The dual-mixed formulation of problem
(\ref{eq:vel_perturbed}) is obtained by seeking $(u,p) \in \Hdiv0 \times
L^2(\Omega)$ such that
$$
\left\{
\begin{aligned}
& \int_\Omega{\nabla \cdot u \delta p \ dx} + \int_\Omega{\lambda p \delta p \
dx} =
\int_\Omega{\Fout \delta p \ dx} \\
& \int_\Omega{u \cdot \delta u \ dx} - \int_\Omega{p \nabla \cdot \delta u \
dx} = 0
\end{aligned}\right.
\quad , \quad \forall (\delta u,\delta p) \in \Hdiv0 \times L^2(\Omega) .
$$
Hence, the bilinear and linear forms associated with the variational
formulation of the problem respectively read as
\begin{gather}
a_{\textrm{vel}}(\{u,p\},\{\delta u,\delta p\}) = \int_\Omega{u \cdot \delta u
\ dx} -
\int_\Omega{p \nabla \cdot \delta u
\ dx} - \int_\Omega{\nabla \cdot u \delta p \ dx} - \int_\Omega{\lambda p
\delta p \ dx}
\label{eq:bilinear_vel_classic} \\
l_{\textrm{vel}}(\{\delta u,\delta p\}) = -\int_\Omega{\Fout \delta p \ dx}
\label{eq:linear_vel_classic}
\end{gather}
To overcome the constraint due to the LBB compatibility condition that the
Finite Element spaces have to fulfill \cite{BoffiBrezziFortin}, several
stabilization approaches have been proposed in the literature over the years
and in this work we consider a strategy inspired by the Galerkin Least-Squares
method and known as CGLS \cite{Correa20081525}.

\subsubsection{Galerkin Least-Squares stabilization}

The GLS formulation relies on adding one or more quantities to the
bilinear form of the problem under analysis in order for the resulting
bilinear form to be strongly consistent and stable.
Let $L$ be the abstract operator for the Boundary Value Problem $L\varphi = g$.
We introduce the solution $\varphi_h$ of the corresponding problem
discretized via the Finite Element Method. The GLS stabilization term reads as
$$
\Lgls(\varphi_h,g;\psi_h) = d \int_\Omega{(L\varphi_h - g)L\psi_h \
dx}.
$$

\subsubsection*{GLS formulation of Darcy's law}

\noindent Following the aforementioned framework, we have
$$
L_1(\{u,p\}) = g_1 \qquad \text{with} \qquad L_1(\{u,p\}) \coloneqq u + \nabla p ,
\qquad g_1 \coloneqq 0.
$$
Thus the GLS term associated with Darcy's law reads as
\begin{equation}
\Lgls_1(\{u_h,p_h\},g_1;\{\delta u_h,\delta p_h\}) = d_1 \int_\Omega{(
u_h
+
\nabla p_h)
\cdot (\delta u_h + \nabla \delta p_h) \ dx}.
\label{eq:GLS_darcy}
\end{equation}

\subsubsection*{GLS formulation of the mass balance equation}

\noindent The equation describing the mass equilibrium may be rewritten as
$$
L_2(\{u,p\}) = g_2 \qquad \text{with} \qquad L_2(\{u,p\}) \coloneqq \nabla
\cdot u + \lambda p, \qquad g_2 \coloneqq \Fout.
$$
Consequently, the Least-Squares stabilization term has the following form
\begin{gather}
\hspace{-12pt} \Lgls_2(\{u_h,p_h\},g_2;\{\delta u_h,\delta p_h\}) =
\underbrace{d_2
\int_\Omega{(\nabla
\cdot u_h + \lambda p_h) (\nabla \cdot \delta u_h + \lambda \delta p_h) \ dx}}
-
\underbrace{d_2
\int_\Omega{\Fout(\nabla \cdot \delta u_h + \lambda \delta p_h) \ dx}}.
\label{eq:GLS_mass} \\
\hspace{130pt} \Lgls_{2a}(\{u_h,p_h\},g_2;\{\delta u_h,\delta p_h\})
\hspace{30pt}
\Lgls_{2l}(\{u_h,p_h\},g_2;\{\delta u_h,\delta p_h\})
\notag
\end{gather}

\subsubsection*{GLS formulation of the curl of Darcy's law}

\noindent Let us consider the rotational component of Darcy's law. Since $p$ is
a scalar
field, $\nabla \times (\nabla p) = 0$ and we get
$$
L_3(\{u,p\}) = g_3 \qquad \text{with} \qquad L_3(\{u,p\}) \coloneqq \nabla \times u, \qquad g_3 \coloneqq 0.
$$
Thus, the GLS term associated with the curl of Darcy's law reads as
\begin{equation}
\Lgls_3(\{u_h,p_h\},g_3;\{\delta u_h,\delta p_h\}) = d_3
\int_\Omega{(\nabla
\times u_h)
(\nabla \times \delta u_h) \ dx}.
\label{eq:GLS_curl}
\end{equation}

\subsubsection*{The stabilized CGLS dual-mixed formulation}

\noindent The stabilized CGLS formulation arises by combining the previous
terms. In
particular, we consider the bilinear form (\ref{eq:bilinear_vel_classic}), we
subtract the Least-Squares stabilization (\ref{eq:GLS_darcy}) for Darcy's law
and we add the corresponding GLS terms (\ref{eq:GLS_mass}) and
(\ref{eq:GLS_curl}) for the mass balance equation and the curl of Darcy's law
itself. In a similar fashion, we assemble the linear form for the stabilized
problem, starting from (\ref{eq:linear_vel_classic}).
The resulting CGLS formulation of problem (\ref{eq:vel_perturbed}) has the
following form:
\begin{gather}
\begin{aligned}
a_{\textrm{CGLS}}(\{u_h,p_h\},\{\delta u_h,\delta p_h\}) = &
a_{\textrm{vel}}(\{u_h,p_h\},\{\delta u_h,\delta p_h\}) -
\Lgls_1(\{u_h,p_h\},g_1;\{\delta u_h,\delta p_h\}) \\
& + \Lgls_{2a}(\{u_h,p_h\},g_2;\{\delta u_h,\delta p_h\}) +
\Lgls_3(\{u_h,p_h\},g_3;\{\delta u_h,\delta p_h\})
\end{aligned}
\label{eq:bilinear_vel_cgls} \\
l_{\textrm{CGLS}}(\{\delta u_h,\delta p_h\}) = l_{\textrm{vel}}(\{\delta
u_h,\delta p_h\}) +
\Lgls_{2l}(\{u_h,p_h\},g_2;\{\delta u_h,\delta p_h\})
\label{eq:linear_vel_cgls}
\end{gather}
To accurately approximate problem
(\ref{eq:bilinear_vel_cgls})-(\ref{eq:linear_vel_cgls}), we consider the
product
space $RT_0 \times \mathbb{P}_0$, that is we use lowest-order Raviart-Thomas
Finite Element for the velocity and piecewise constant functions for the
pressure.

\subsection{Streamline Upwind Petrov Galerkin for the dynamics of gases and
liquid water}

The numerical approximation of pure advection and advection-diffusion
transient problems has to be carefully handled in order to retrieve accurate
solutions. It is well-known in the literature \cite{Bochev20042301} that
classical Finite Element Method suffers from poor accuracy when dealing with
steady-state advection and advection-diffusion problems and requires the
introduction of numerical stabilization to construct a strongly consistent
scheme.
When moving to transient advection and advection-diffusion problems, time-space
elements are the most natural setting to develop stabilized methods
\cite{JOHNSON1984285}.

Let $L_{\textrm{ad}}$ be the abstract operator to model an advection -
respectively advection-diffusion - phenomenon. The resulting transient Boundary
Value
Problem may be written as
\begin{equation}
\phi \partial_t \varphi + L_{\textrm{ad}}\varphi = g_{\textrm{ad}} .
\label{eq:abstractAD}
\end{equation}
We consider the variational formulation of (\ref{eq:abstractAD}) by introducing
the corresponding abstract bilinear form $B_{\textrm{ad}}(\varphi,\psi)$ which
will be detailed in next subsections.
Let $\varphi_h$ be the solution of the discretized PDE via the Finite
Element Method. The SUPG stabilization term for the transient problem reads as
$$
\Lsupg(\varphi_h,g_{\textrm{ad}};\psi_h) = d \int_\Omega{(\phi \partial_t
\varphi_h + L_{\textrm{ad}}\varphi_h - g_{\textrm{ad}})
L_{\textrm{ad}}^{\textrm{SS}}\psi_h \ dx}
$$
where $L_{\textrm{ad}}^{\textrm{SS}}$ is the skew-symmetric part of the
advection - respectively advection-diffusion - operator and $d$ is a
stabilization parameter constant in space and time. Let $\mathcal{T}_h$ be the
computational mesh that approximates the domain $\Omega$ and $\ell_K$ the
diameter of each element $K \in \mathcal{T}_h$. We choose:
\begin{equation}
d \coloneqq \frac{1}{2 \| u \|_2} \max_{K \in \mathcal{T}_h} \ell_K .
\label{eq:stabilizationParameter}
\end{equation}
Let $V_h \coloneqq \{ \psi_h \in \mathcal{C}(\overline{\Omega}) \ \text{s.t.} \
\psi_h|_K \in \mathbb{P}_k(K) \ \forall K \in \mathcal{T}_h \}$ be the space of
Lagrangian Finite Element of degree $k$, that is the piecewise polynomial
functions of degree $k$ on each element $K$ of the mesh $\mathcal{T}_h$.
The stabilized SUPG formulation of the advection - respectively
advection-diffusion - problem (\ref{eq:abstractAD}) reads as follows: for all
$t \in (0,\Sfin]$ we seek $\varphi_h(t) \in V_h$ such that
\begin{equation}
\int_\Omega \phi \partial_t \varphi_h \psi_h \ dx +
B_{\textrm{ad}}(\varphi_h,\psi_h) + \Lsupg(\varphi_h,g_{\textrm{ad}};\psi_h) =
\int_\Omega g_{\textrm{ad}} \psi_h \ dx \qquad \forall \psi_h \in V_h .
\label{eq:SUPG}
\end{equation}

Concerning time discretization, it is well-known in the literature that 
implicit schemes tend to increase the overall computational cost associated 
with the solution of a PDE.
Nevertheless, precise approximations of advection and advection-diffusion 
equations via explicit schemes usually require high-order methods and are 
subject to stability conditions that may be responsible of making computation 
unfeasible.
On the contrary, good stability and convergence properties of implicit 
strategies make them an extremely viable option when dealing with complex - 
possibly coupled - phenomena and with equations featuring noisy parameters.
In particular, owing to the coupling of several PDE's, the solution of the 
advection and advection-diffusion equations presented in our model 
turned out to be extremely sensitive to the choice of the involved parameters. 
Being the tuning of the unknown coefficients of the equations one of the main 
goal of the SiViBiR++ project, a numerical scheme unconditionally stable and 
robust to the choice of the discretization parameters is sought.
Within this framework, we consider an implicit Euler scheme for the 
time discretization and we stick to low-order Lagrangian Finite Element functions 
for the space discretization. 
The numerical scheme arising from the solution of equation (\ref{eq:SUPG}) by
means of the aforementioned approximation is known to be stable and to
converge quasi-optimally \cite{Burman20101114}.

In the following subsections, we provide some details on the bilinear and
linear forms involved in the discretization of the advection and
advection-diffusion equations as well as on the formulation of the
associated stabilization terms.

\subsubsection{The case of gases}

Let us consider a generic gas $G$ whose Finite Element counterpart is
named $G_h$. Within the previously introduced framework, we get
$$
L_G G \coloneqq u \cdot \nabla G, \qquad g_G \coloneqq F^G,
\qquad B_G(G,\delta G) = \int_\Omega (u \cdot \nabla G ) \delta G \ dx.
$$
Hence, the SUPG stabilization term reads as
$$
\Lsupg_G(G_h,g_G;\delta G_h) = d \int_\Omega{(\phi \partial_t G_h + u \cdot
\nabla G_h - F^G) (u \cdot \nabla \delta G_h) \ dx} .
$$
By introducing an implicit Euler scheme to approximate the time derivative in
(\ref{eq:SUPG}), we obtain the fully discretized advection problem for the
gas $G$: at $t=t_n$ we seek $G_{h,n+1} \in V_h$ such that $a_G(G_{h,n+1},\delta
G_h) = l_G(\delta G_h) \ \forall \delta G_h \in V_h$, where
\begin{gather}
a_G(G_{h,n+1},\delta G_h) = \int_\Omega \phi G_{h,n+1} (\delta G_h + d \ u
\cdot \nabla \delta G_h) \, dx +
\Delta t \int_\Omega  u \cdot \nabla G_{h,n+1} (\delta G_h + d \ u \cdot
\nabla \delta G_h) \, dx \, ,
\label{eq:variational_a_G} \\
l_G(\delta G_h) = \int_\Omega \phi G_{h,n} (\delta G_h + d \
u \cdot \nabla \delta G_h) \, dx + \Delta t \int_\Omega F^G (\delta G_h + d \
u \cdot \nabla \delta G_h) \, dx .
\label{eq:variational_l_G}
\end{gather}
\begin{rmrk}
As the authors highlight in \cite{Bochev20042301}, from a practical point
of view the implementation of
(\ref{eq:variational_a_G})-(\ref{eq:variational_l_G}) may not be
straightforward due to the non-symmetric mass matrix resulting from the
discretization of the first term in (\ref{eq:variational_a_G}).
\end{rmrk}
\noindent By considering $F^O = F^N = 0$, we get the linear forms associated
with the dynamic of oxygen and nitrogen:
\begin{equation}
\begin{split}
& l_O(\delta O_h) = \int_\Omega \phi O_{h,n} (\delta O_h + d \
u \cdot \nabla \delta O_h) \, dx \, , \\
& l_N(\delta N_h) = \int_\Omega \phi N_{h,n} (\delta N_h + d \
u \cdot \nabla \delta N_h) \, dx .
\end{split}
\label{eq:l_OxyNitro}
\end{equation}
In a similar way, when $F^j = -c_j \partial_t \Corg \, \ j=M,C$, we obtain the
linear forms for methane and carbon dioxide:
\begin{equation}
\begin{split}
& l_M(\delta M_h) = \int_\Omega \phi M_{h,n} (\delta M_h + d \
u \cdot \nabla \delta M_h) \, dx - \int_\Omega c_M
(\Corg_{h,n+1}-\Corg_{h,n}) (\delta M_h + d \
u \cdot \nabla \delta M_h) \, dx \, , \\
& l_{\Cdx}(\delta C_h) = \int_\Omega \phi \Cdx_{h,n} (\delta C_h + d \
u \cdot \nabla \delta C_h) \, dx - \int_\Omega c_C
(\Corg_{h,n+1}-\Corg_{h,n}) (\delta C_h + d \
u \cdot \nabla \delta C_h) \, dx .
\end{split}
\label{eq:l_MCdx}
\end{equation}
Eventually, the dynamic of the water vapor is obtained when considering $F^h =
-c_h \partial_t \Corg - \Fcond$:
\begin{equation}
l_h(\delta h_h) = \int_\Omega \phi h_{h,n} (\delta h_h + d \
u \cdot \nabla \delta h_h) \, dx - \int_\Omega c_h
(\Corg_{h,n+1}-\Corg_{h,n}) (\delta h_h + d \
u \cdot \nabla \delta h_h) \, dx  - \Delta t \int_\Omega \Fcond (\delta h_h + d
\ u \cdot \nabla \delta h_h) \, dx .
\label{eq:l_vapor}
\end{equation}

\subsubsection{The case of liquid water}

The dynamic of the liquid water being described by an advection-diffusion
equation, the SUPG framework may be written in the following form:
$$
L_w w \coloneqq u_w \cdot \nabla w - k_w \Delta w, \qquad g_w
\coloneqq \Fcond + \Fin,
\qquad B_w(w,\delta w) = \int_\Omega \Big( (u_w \cdot \nabla w ) \delta w + k_w
\nabla w \cdot \nabla \delta w \Big) dx.
$$
Hence the stabilization term reads as
$$
\Lsupg_w(w_h,g_w;\delta w_h) = d_w \int_\Omega{(\phi \partial_t w_h + u_w \cdot
\nabla w_h - k_w \Delta w_h - \Fcond - \Fin) (u_w \cdot \nabla \delta w_h) \
dx}
$$
where $d_w$ is obtained by substituting $u_w$ in
(\ref{eq:stabilizationParameter}). We obtain the fully discretized problem in
which at each $t=t_n$ we seek $w_{h,n+1} \in V_h$ such that
\begin{gather}
a_w(w_{h,n+1},\delta w_h) = l_w(\delta w_h) \qquad \forall \delta w_h \in V_h
\notag \\
\begin{aligned}
a_w(w_{h,n+1},\delta w_h) = & \int_\Omega \phi w_{h,n+1} (\delta w_h + d_w
\ u_w \cdot \nabla \delta w_h) \, dx +
\Delta t \int_\Omega  u_w \cdot \nabla w_{h,n+1} (\delta w_h + d_w \ u_w \cdot
\nabla \delta w_h) \, dx \\
& + \Delta t \int_\Omega  k_w \nabla w_{h,n+1} \cdot \nabla \delta w_h \, dx -
\Delta t \int_\Omega  k_w \Delta w_{h,n+1} (d_w \ u_w \cdot \nabla \delta w_h)
\, dx \, ,
\end{aligned}
\label{eq:variational_a_w} \\
l_w(\delta w_h) = \int_\Omega \phi w_{h,n} (\delta w_h + d_w \
u_w \cdot \nabla \delta w_h) \, dx +
\Delta t \int_\Omega (\Fcond + \Fin) (\delta w_h + d_w \ u_w \cdot \nabla
\delta
w_h) \, dx .
\label{eq:l_liquid}
\end{gather}

\section{Numerical results}
\label{ref:simulations}

In this section we present some preliminary numerical simulations to test the
proposed model. The SiViBiR++ project implements the discussed numerical methods for 
the approximation of the phenomena inside a bioreactor landfill. It is based on 
a C++ library named \Feel which provides a framework to solve PDE's using the 
Finite Element Method \cite{Feel2012}.

\subsection{\Feel}

\Feel stands for \emph{Finite Element Embedded Language in C++} and
is a C++ library for the solution of Partial Differential Equations using
generalized Galerkin methods.
It provides a framework for the implementation of advanced numerical methods to
solve complex systems of PDE's.
The main advantage of \Feel for the applied mathematicians and engineers
community relies on its design based on the Domain Specific Embedded Language
(DSEL) approach \cite{Prud'homme:2006:DSE:1376891.1376895}.
This strategy allows to decouple the difficulties encountered by the scientific
community when dealing with libraries for scientific computing.
As a matter of fact, DSEL provides a high-level language to handle mathematical
methods without loosing abstraction.
At the same time, due to the continuing evolution of the state-of-the-art
techniques in computer science (e.g. new standards in programming languages,
parallel architectures, etc.) the choice of the proper tools in scientific
computing may prove very difficult.
This is even more critical for scientists who are not specialists in computer science
and have to reach a compromise between user-friendly interfaces and high performances.

\Feel proposes a solution to hide these difficulties behind a
user-friendly language featuring a syntax that mimics the mathematical formulation
by using a a much more common low-level language, namely C++. \\
Moreover, \Feel integrates the latest C++ standard - currently C++14 - 
and provides seamless parallel tools to handle mathematical operations such
as projection, integration or interpolation through C++ keywords.
\Feel is regularly tested on High Performance Computing facilities such as the
PRACE research infrastructures (e.g. Tier-0 supercomputer CURIE, Supermuc, etc.)
via multidisciplinary projects mainly gathered in the \Feel consortium.

The embedded language provided by the \Feel framework represents a powerful engineering
tool to rapidly develop and deploy production of ready-to-use softwares as well as
prototypes.
This results in the possibility to treat physical and engineering applications
featuring complex coupled systems from early-stage exploratory analysis till the
most advanced investigations on cutting-edge optimization topics.
Within this framework, the use of \Feel for the simulation of the dynamic
inside a bioreactor landfill seemed promising considering the complexity of
the problem under analysis featuring multiphysics phenomena at different space
and time scales.

A key aspect in the use of \Feel for industrial applications like the one
presented in this paper is the possibility to operate on parallel
infrastructures without directly managing the MPI communications.
Here, we
briefly recall the main steps for the parallel simulation of the dynamic
inside a bioreactor landfill highlighting the tools involved in \Feel and in
the external libraries linked to it:
\begin{itemize}
 \item we start by constructing a computational mesh using \Gmsh \cite{GMSH};
 \item a mesh partition is generated using Chaco or Metis and
additional information about ghost cells with communication between neighbors
is provided \cite{GMSH};
 \item \Feel generates the required parallel data structures and create a
table with global and local views for the Degrees of Freedom;
 \item \Feel assembles the system of PDE's starting from the variational
formulations and the chosen Finite Element spaces;
 \item the algebraic problem is solved using the efficient solvers and
preconditioners provided by PETSc \cite{petsc2015}.
\end{itemize}
A detailed description of the high performance framework within \Feel is
available in \cite{ChabanneThesis}.

\subsection{Geometric data}\label{ssec:geometric_data}

\begin{figure}[tbp]
\centering
\subfloat[Geometry of the alveolus.]
{
    \includegraphics[width=0.31\columnwidth]{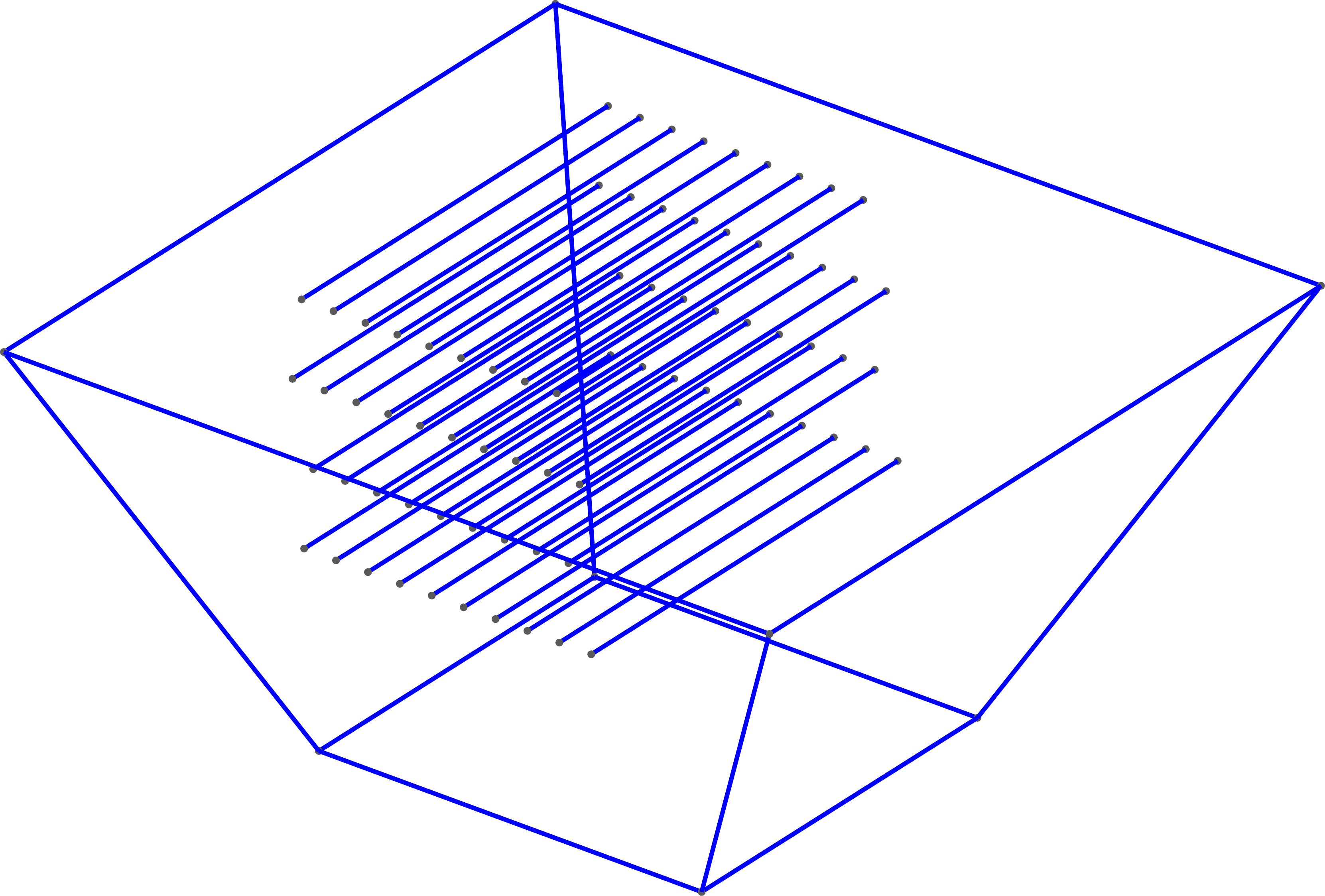}
\label{fig:domain7b}
}
\hfill
\subfloat[Surface mesh.]
{
    \includegraphics[width=0.31\columnwidth]{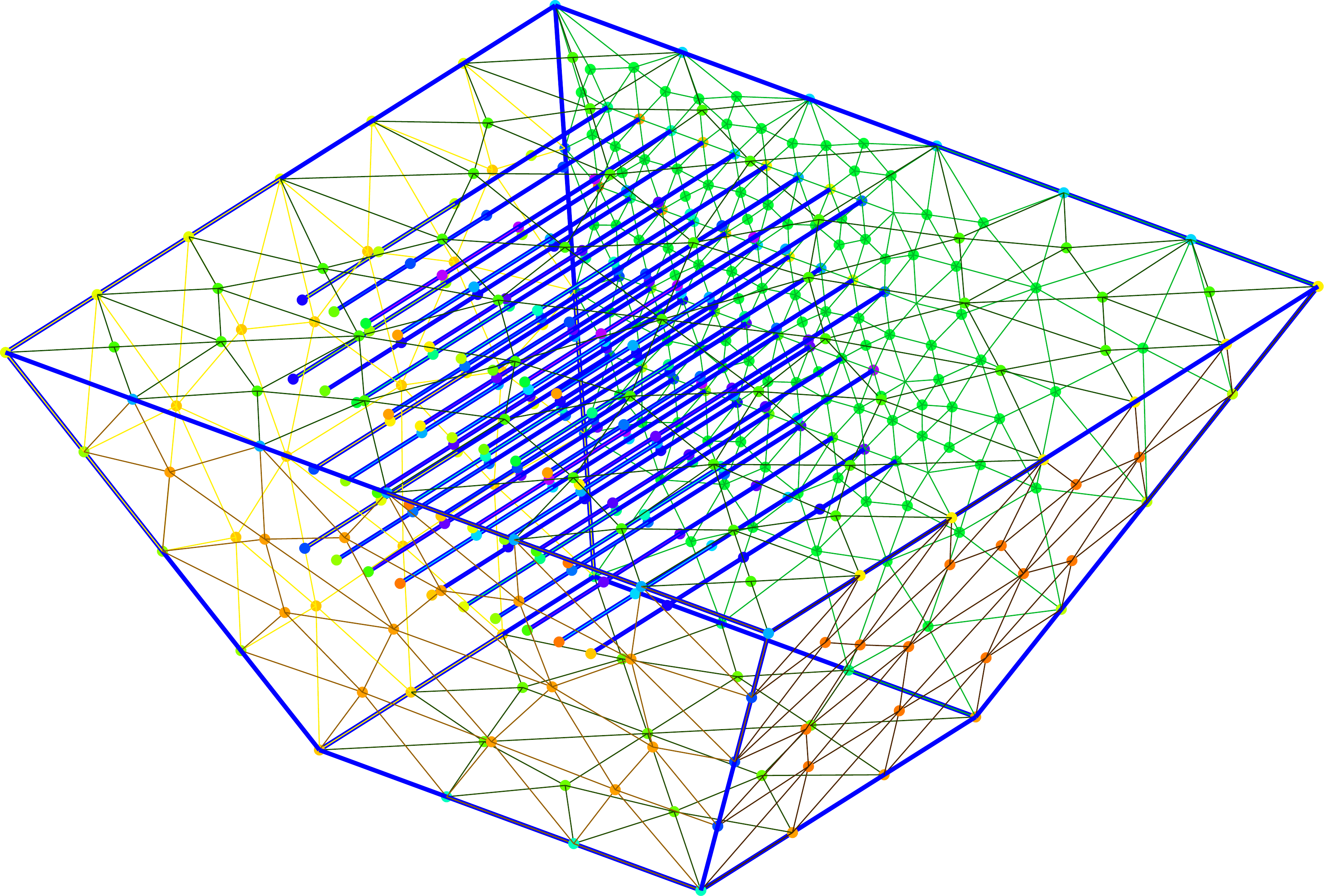}
\label{fig:mesh2d7b}
}
\hfill
\subfloat[Volume mesh.]
{
    \includegraphics[width=0.31\columnwidth]{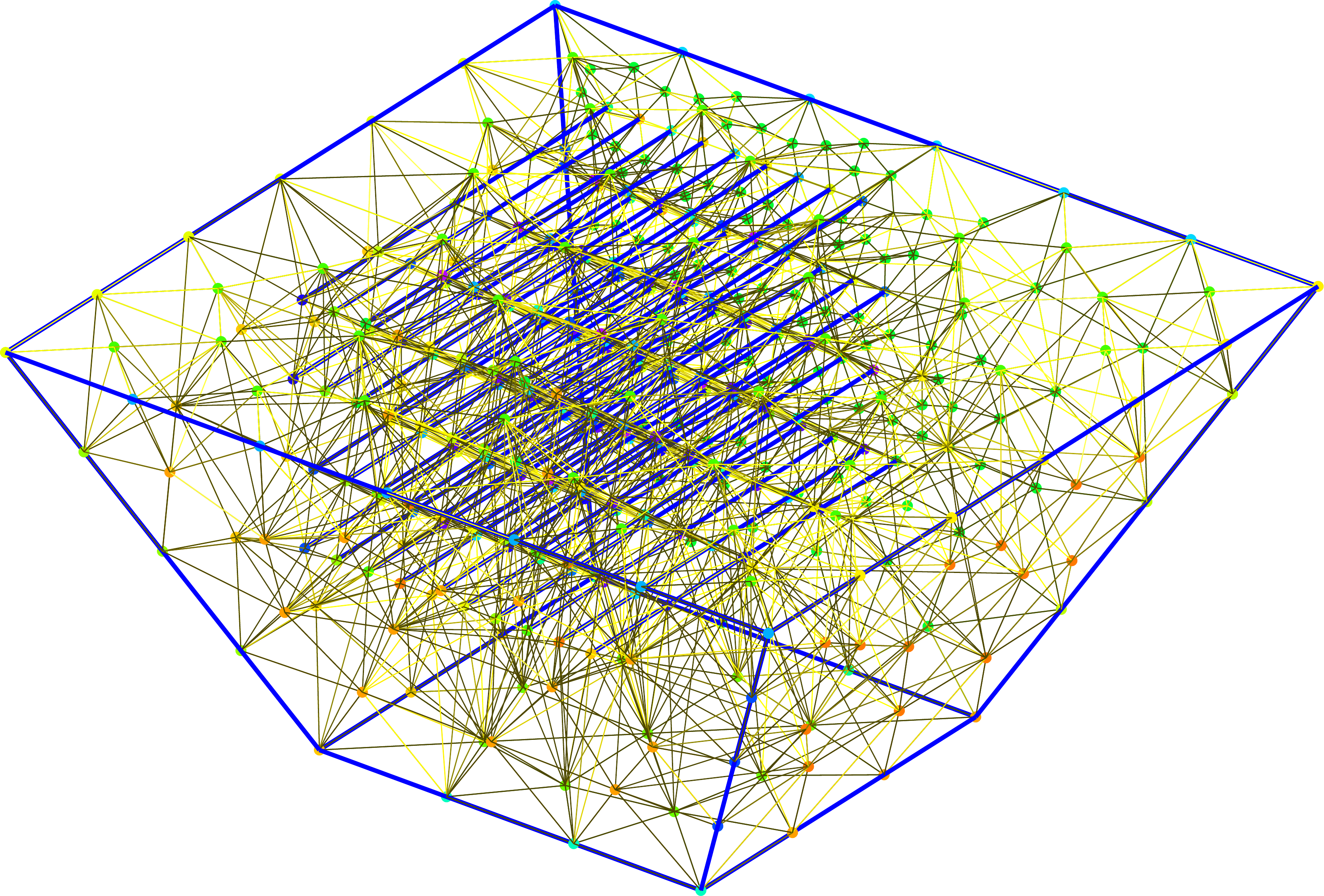}
\label{fig:mesh3d7b}
}
\caption{Geometry of the alveolus and computational mesh generated using \Gmsh. Each layer of 
1D lines alternatively represents a set of water injectors or a group of extraction drains for 
the gas.}
\label{fig:geomesh}
\end{figure}

We consider a reverse truncated pyramid domain as in figure \ref{fig:pyramid}.
The base of the domain measures $\SI{90}{m} \times \SI{90}{m}$ and its height
is $\SI{90}{m}$. All the lateral walls feature a slope of $\pi/6$. The alveolus
counts $20$ extraction drains organized on $2$ levels and $20$ injection pipes,
distributed on $2$ levels as well. All the pipes are $\SI{25}{m}$ long and are
modeled as 1D lines since their diameters are of order $\SI{e-1}{m}$. A
simplified scheme of the alveolus under analysis is provided in figure
\ref{fig:domain7b} and the corresponding computational domain is obtained by
constructing a triangulation of mesh size $\SI{5}{m}$ (Fig. \ref{fig:mesh3d7b}).

\subsection{Heuristic evaluation of the unknown constants in the model}
\label{ref:heuristic}

The model presented in this paper features a large set of unknown variables
(i.e. diffusion coefficients, scaling factors, ...) whose role is crucial to
obtain realistic simulations. In this section, we propose a first set of values
for these parameters that have been heuristically deduced by means of some
qualitative and numerical considerations. A major improvement of the model is
expected by a more rigorous tuning of these parameters which will be
investigated in future works.

\subsubsection*{Porous medium}
\noindent The physical and chemical properties of the
bioreactor considered as a porous environment have been derived by experimental
results in the literature. In particular, we consider a porosity $\Phi =
\SI{0.3}{}$ and a permeability $D=\SI{e-11}{m^2}$.
\subsubsection*{Bacteria and organic carbon}
\noindent We consider both the concentration of bacteria $b$ and of organic
carbon $\Corg$ as non-dimensional quantities in order to estimate their
evolution. Thus we set $b_0=1$ and $\Corg_0=1$ and we derive the values
$a_b=\SI{e-5}{m^6\kg^{-2}\day^{-1}}$ and $c_b=1$ respectively for the rate of
consumption of the organic carbon and for the rate of reproduction of bacteria.
Within this framework and under the optimal conditions of reaction prescribed
by (\ref{eq:psi1psi2}), $\Corg$ decreases to $2\permil$ of its initial value
during the lifetime of the bioreactor whereas bacteria concentration $b$ remains
bounded ($b < 2b_0$).
\subsubsection*{Temperature}
\noindent As per experimental data, the optimal temperature for the methanogenic
fermentation to take place is $T_{opt} = \SI{35}{\celsius}=\SI{308}{K}$
with an admissible variation  of temperature of
$\pm A_T=\SI{20}{\celsius}=\SI{20}{K}$ to guarantee the survival of bacteria.
The factor $c_T$ represents the heat produced per unit of consumed organic
carbon and per unit of time and is estimated to $c_T =
\SI{e2}{K.m^2\day^{-1}}$. The thermal conductivity of the waste inside the
alveolus is fixed at $k_T = \SI{9e-2}{m^{2}\day^{-1}}$. In order to impose
realistic boundary conditions for the heat equation, we consider different
values for the temperature of the lateral surface of the alveolus
$T_g=\SI{5}{\celsius}=\SI{278}{K}$ and the one of the top geomembrane
$T_m=\SI{20}{\celsius}=\SI{293}{K}$.
\subsubsection*{Water}
\noindent The production of methane takes place only when less
than $10\%$ of water is present inside the bioreactor. Since the alveolus is
completely flooded when $w =\SI{1000}{\kg.m^{-3}}$, we get that $w_\text{max}
= \SI{100}{\kg.m^{-3}}$. The vertical velocity drift due to gravity is
estimated from Darcy's law to the value $\|u_w\| = \SI{2,1}{m.\day^{-1}}$ and
the diffusion coefficient is set to $k_w = \SI{8.6e-2}{m^2\day^{-1}}$.
\subsubsection*{Phase transition}
\noindent In order to model the phase transitions, we have to consider the
critical values of the pressure associated with evaporation and condensation.
The Rankine law for the vapor pressure of water is approximated using the
following values: $P_0 = \SI{133.322}{\Pa}$, $s_0 = 20.386$ and $s_1 = 5132 K$
and its linearization arises when considering $H_0 = \SI{-9.56e4}{\Pa}$ and $H_1
= \SI{337.89}{\Pa.K^{-1}}$ for the range of temperatures
$[\SI{288}{K};\SI{328}{K}]$. Moreover, we set the value $c_{h\to w}$ that
represents the speed for the condensation of vapor to liquid water: $c_{h\to w}
= \SI{e-1}{\day^{-1}}$.
\subsubsection*{Gases}
\noindent We consider a gas mixture made of methane, carbon dioxide and water
vapor. Its dynamical viscosity is set to
$\mu_\text{gas}=\SI{1,3}{Pa.\day^{-1}}$. Other parameters involved represent the
rate of production of a specific gas (methane, carbon dioxide and water vapor)
per unit of consumed organic carbon and per unit of time: $c_M =
\SI{1.8e7}{\kg.m^{-3}}$; $c_C = \SI{2.6e7}{\kg.m^{-3}}$; $c_h =
\SI{2,5e6}{\kg.m^{-3}}$.

\begin{rmrk}
A key aspect in the modeling of a bioreactor landfill is the possibility to
adapt the incoming flow of water and leachates $\Jin$ and the outgoing flow of
biogas $\Jout$. These values are user-defined parameters which are kept
constant to $\SI{258}{m^3.d^{-1}}$ for the simulations in this paper but should
act as control variables in the framework of the forecast and optimization
procedures described in the introduction.
\end{rmrk}

\subsection{A preliminary test case}

In this section we present some preliminary numerical results obtained by using the 
SiViBiR++ module developed in \Feel to solve the equations presented in section 
\ref{ref:math_model} using the numerical schemes discussed in section \ref{ref:numerical_model}.
In particular, we remark that in all the following simulations we neglect the 
effects due to the gas and fluid dynamics inside the bioreactor landfill.
Though this choice limits the ability of the discussed results to correctly 
describe the complete physical behavior of the system, this simplification is a necessary 
starting point for the validation of the mathematical model in section \ref{ref:math_model}.
As a matter of fact, the equations describing the gas and fluid dynamics feature 
several unknown parameters whose tuning - independent and coupled with one another - 
has to be accurately performed before linking them to the problems modeling the 
consumption of organic carbon and the evolution of the temperature. \\
Thus, here we restrict our numerical simulations to two main phenomena occurring 
inside the bioreactor landfill: first, we describe the consumption of organic 
carbon under some fixed optimal conditions of humidity and temperature; then we 
introduce the evolution in time of the temperature and we discuss the behavior 
of the coupled system given by equations 
(\ref{eq:organic_carbon_eq})-(\ref{eq:heat_eq}).

The test cases are studied in the computational domain introduced 
in section \ref{ssec:geometric_data}: in particular, we consider the triangulation 
of mesh size $\SI{5}{m}$ in figure \ref{fig:mesh3d7b} and we set the unit 
measure for the time evolution to $\Delta t=\SI{365}{\day}$. The final time 
for the simulation is $\Sfin=40$ years.
The parameters inside the equations are set according to the values in 
section \ref{ref:heuristic} but a thorough investigation of these quantities has to 
be performed to verify their accuracy.
The computations have been executed using up to $32$ processors and below we 
present some simulations for the aforementioned preliminary test cases.

\subsubsection*{Evolution of the organic carbon under optimal hydration 
and temperature conditions}

First of all, we consider the case of a single uncoupled equation, that is  
the evolution of the organic carbon in a scenario in which the concentration of water 
and the temperature are fixed. 
Starting from equation (\ref{eq:organic_carbon_eq}), we assume fixed optimal conditions for 
the humidity and the temperature inside the bioreactor landfill. 
We set a fixed amount of water $w=\frac{\wMax}{2}$ inside the alveolus and 
a constant temperature $T=\Topt$. Thus, from (\ref{eq:psi1psi2}) we get
$$
\Psi_1(w) \equiv \frac{\wMax}{4}
,\qquad
\Psi_2(T) \equiv 1
$$
and equation (\ref{eq:organic_carbon_eq}) reduces to
\begin{equation}
\begin{cases}
\displaystyle
(1-\phi) \partial_t \Corg(x,t) = - a_b \ \frac{\wMax}{4} [ b_0 + c_b (\Corg_0 - \Corg(x,t)) ] \ \Corg(x,t) \ , \ & \text{in }\Omega \times (0,\Sfin] \\
\Corg(\cdot,0) = \Corg_0 \ , \ & \text{in }\overline{\Omega}
\end{cases}
\label{eq:optimalCorg}
\end{equation}
It is straightforward to observe that equation (\ref{eq:optimalCorg}) only features one 
unknown variable - namely the organic carbon - since the concentration of bacteria $b(x,t)$ 
is an affine function of the concentration of organic carbon itself (cf. equation (\ref{eq:bacteria})). 

\begin{figure}[htbp]
\centering
\subfloat[Lifetime: $1$ year.]
{
\includegraphics[width=0.48\columnwidth]{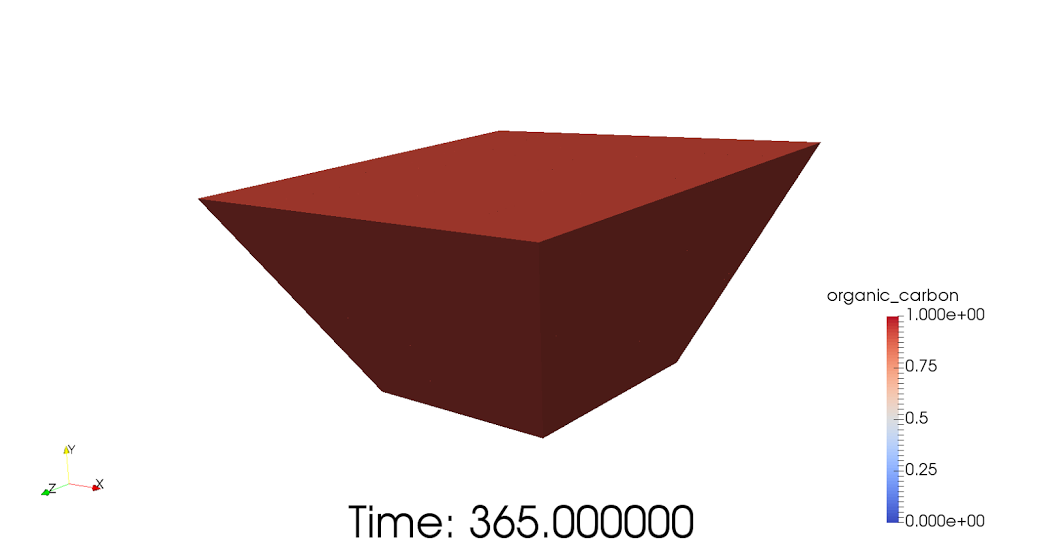}
\label{fig:corg1}
}
\hfill
\subfloat[Lifetime: $10$ years.]
{
\includegraphics[width=0.48\columnwidth]{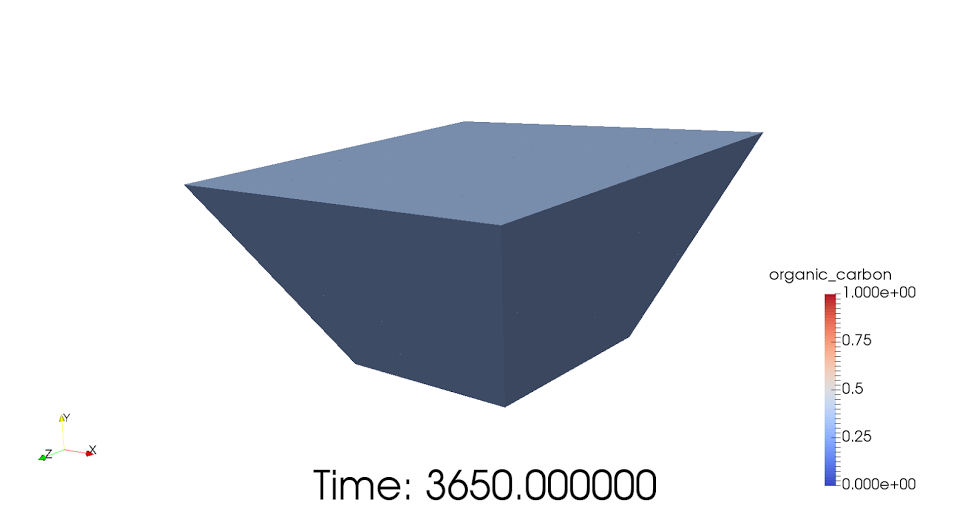}
\label{fig:corg10}
}

\subfloat[Lifetime: $20$ years.]
{
\includegraphics[width=0.48\columnwidth]{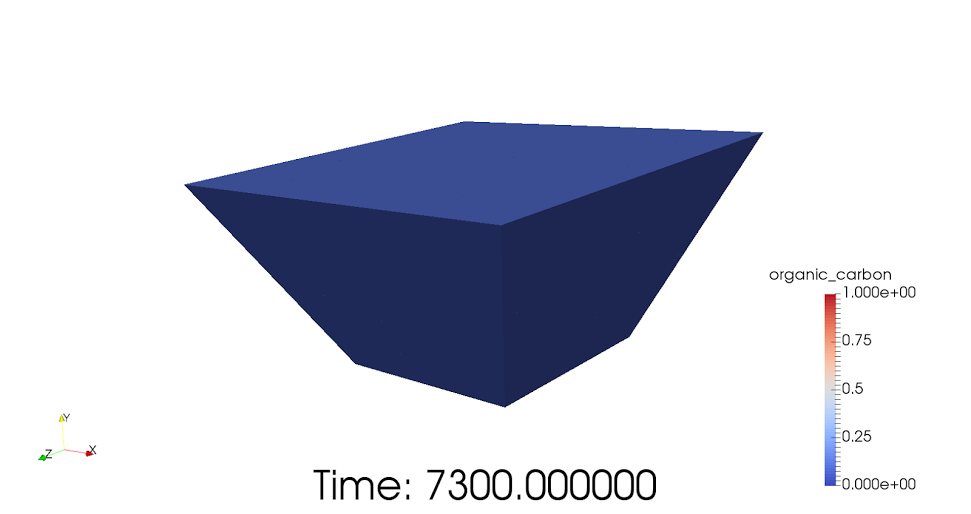}
\label{fig:corg20}
}
\hfill
\subfloat[Lifetime: $40$ years.]
{
\includegraphics[width=0.48\columnwidth]{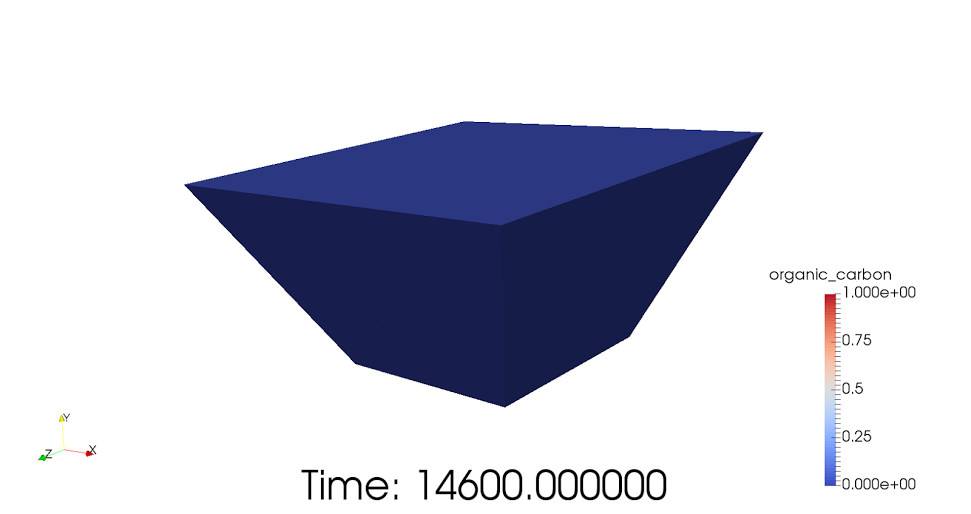}
\label{fig:corg40}
}
\caption{Evolution of the organic carbon inside the alveolus at $t=1$ year, 
$t=10$ years, $t=20$ years and $t=40$ years.}
\label{fig:Corg_snapshots}
\end{figure}

\begin{figure}[hbtp]
\centering
\includegraphics[width=0.5\columnwidth]{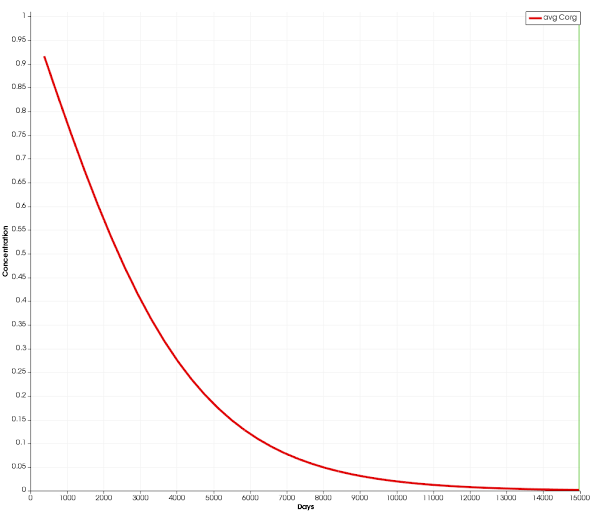}
\caption{Evolution of the organic carbon between $t=0$ and $t=40$ years 
under optimal hydration and temperature conditions.}
\label{fig:Corg_evolution}
\end{figure}

As stated in section \ref{ssec:num_organic_heat_eq}, the key aspect in the solution of equation 
(\ref{eq:optimalCorg}) is the handling of the non-linear reaction term. In order to numerically treat this term 
as described, at each time step we need the value $\Corg_n$ at the previous iteration to compute 
the semi-implicit quantity $\Corg_{n+1}\Corg_n$. 
To provide a suitable value of $\Corg_n$ during the first iteration, we solve a linearization 
of equation (\ref{eq:optimalCorg}) and we use the corresponding solution to evaluate $\Corg_{n+1}\Corg_n$. \\
The initial concentration of organic carbon inside the bioreactor is set to $1$ and in 
figure \ref{fig:Corg_snapshots} we observe several snapshots of the quantity of organic carbon 
inside the alveolus at time $t=1$ year, $t=10$ years, $t=20$ years and $t=40$ years. 
At the end of the life of the facility, the amount of organic carbon inside the alveolus is 
$\Corg=\SI{2.e-3}{}$.
In figure \ref{fig:Corg_evolution}, we plot the evolution of the overall 
quantity of organic carbon with respect to time. 
In particular, as expected we observe that $\Corg$ decreases in time as 
the methanogenic fermentation takes place. 

\subsubsection*{Evolution of the temperature as a function of organic carbon}

\begin{figure}[htbp]
\centering
\subfloat[Lifetime: $1$ year.]
{
\includegraphics[width=0.48\columnwidth]{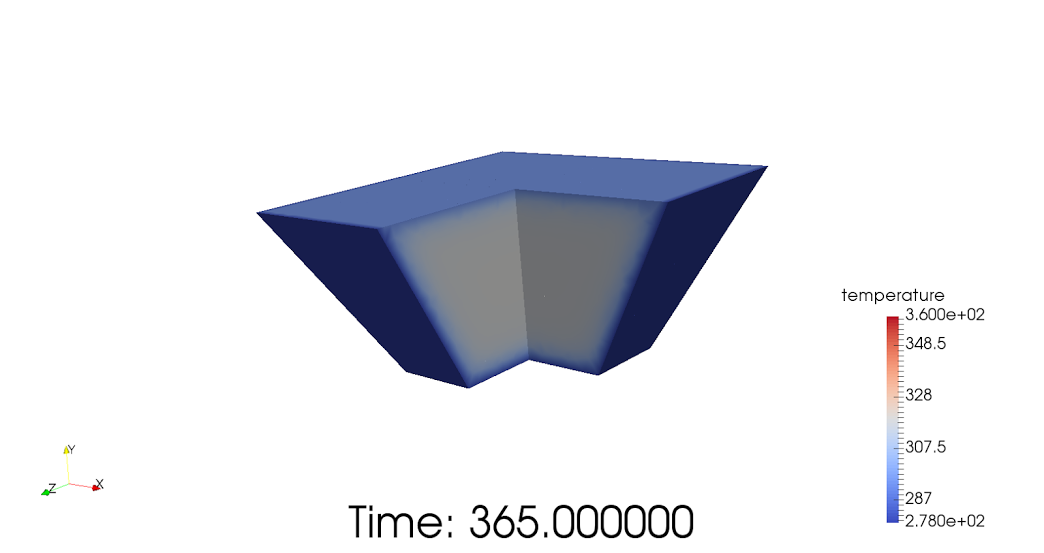}
\label{fig:temp1}
}
\hfill
\subfloat[Lifetime: $10$ years.]
{
\includegraphics[width=0.48\columnwidth]{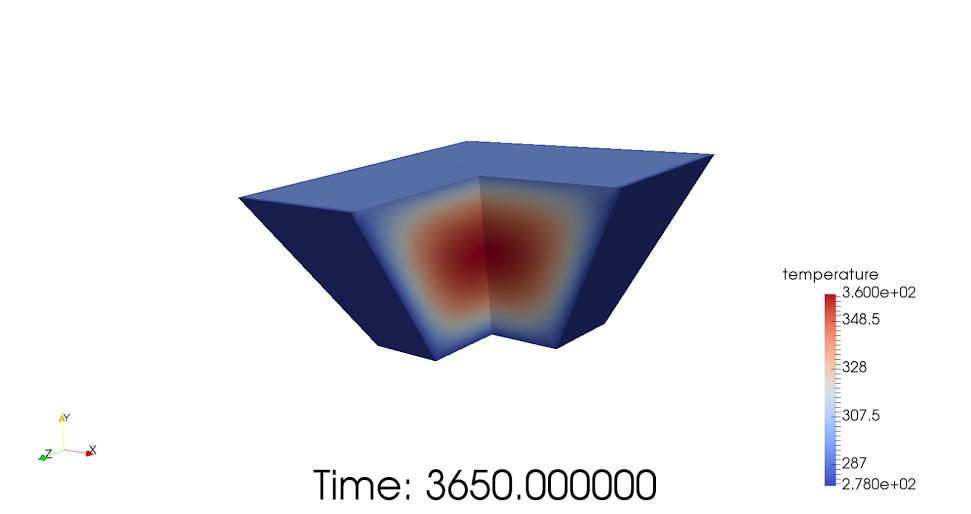}
\label{fig:temp10}
}

\subfloat[Lifetime: $20$ years.]
{
\includegraphics[width=0.48\columnwidth]{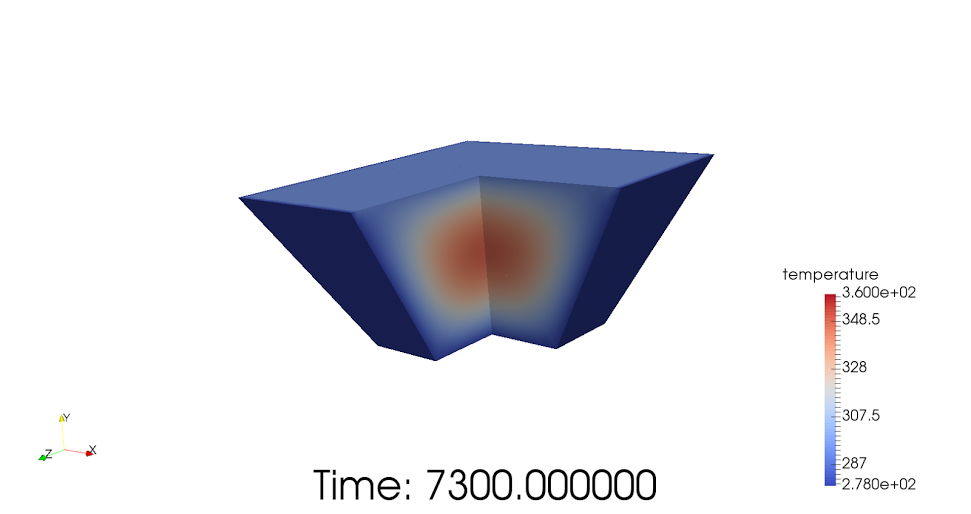}
\label{fig:temp20}
}
\hfill
\subfloat[Lifetime: $40$ years.]
{
\includegraphics[width=0.48\columnwidth]{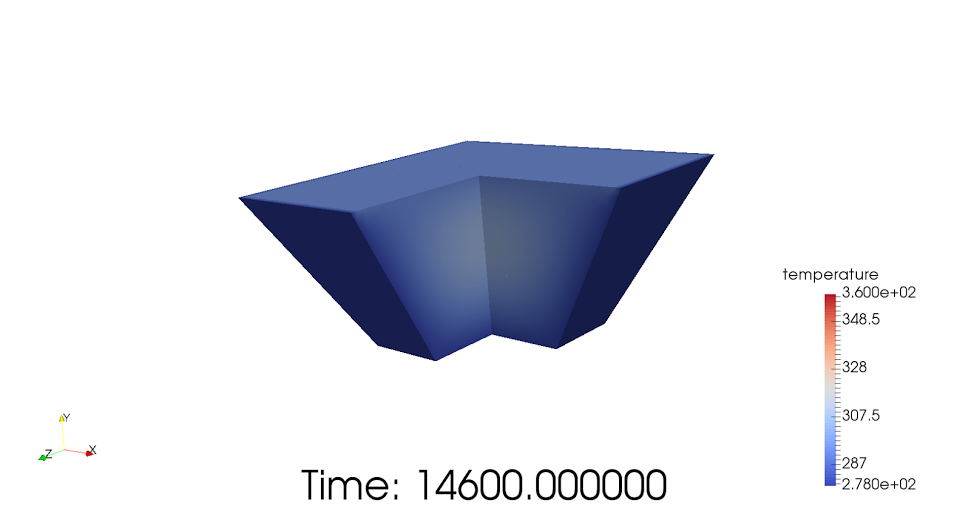}
\label{fig:temp40}
}
\caption{Evolution of the temperature inside the alveolus at $t=1$ year, 
$t=10$ years, $t=20$ years and $t=40$ years.}
\label{fig:temp_alone_snapshots}
\end{figure}

\begin{figure}[hbtp]
\centering
\subfloat[Evolution of the organic carbon.]
{
\includegraphics[width=0.48\columnwidth]{./Fig5-7a.png}
\label{fig:carb_ev}
}
\hfill
\subfloat[Evolution of the temperature.]
{
\includegraphics[width=0.48\columnwidth]{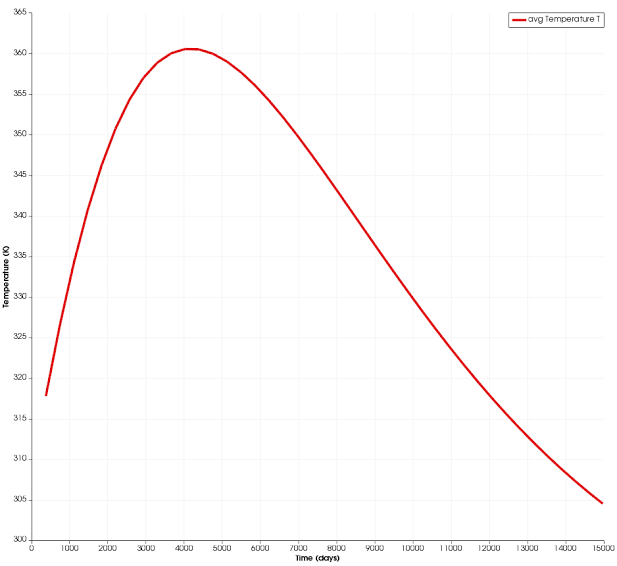}
\label{fig:temp_ev}
}
\caption{Evolution of the quantity of organic carbon and temperature inside the alveolus 
between $t=0$ and $t=40$ years.}
\label{fig:temp_carb}
\end{figure}

We introduce a novel variable depending both on space and time to model the 
temperature inside the bioreactor.
Figure \ref{fig:temp_alone_snapshots} presents the snapshots of the value of the 
temperature in a section of the alveolus under analysis at time $t=1$ year, $t=10$ 
years, $t=20$ years and $t=40$ years. 
These result from the solution of equation (\ref{eq:heat_eq}) for a given 
trend of the organic carbon. 
Let us consider the evolution of the organic carbon obtained from the previous 
test case. The corresponding profile is given by
$$
\Corg(x,t) = e^{-\alpha t} \qquad , \qquad \alpha = \SI{e-3}{}.
$$
We observe that the consumption of 
organic carbon by means of the chemical reaction (\ref{eq:chemical_reaction}) 
is responsible for the generation of heat in the middle of the domain. As expected by 
the physics of the problem, the heat tends to diffuse towards the external boundaries 
where the temperature is lower.
After a first phase which lasts approximately $10$ years in which the methanogenic 
process produces heat and the temperature rises, the consumption of organic carbon 
slows down and the temperature as well starts to decrease until the end-of-life of 
the facility (Fig. \ref{fig:carb_ev}-\ref{fig:temp_ev}).

\subsubsection*{Evolution of the coupled system of organic carbon and 
temperature under optimal hydration condition}

Starting from the previously discussed cases, we now proceed to the coupling of 
the organic carbon with the temperature.
We keep the optimal hydration condition as in the previous simulations - 
that is $w=\frac{\wMax}{2}$ - and we consider the solution of the coupled equations 
(\ref{eq:organic_carbon_eq})-(\ref{eq:heat_eq}).

From the numerical point of view, this scenario introduces several difficulties, 
mainly due to the fact that the two equations are now dependent on one another. 
As mentioned in section \ref{ref:numerical_model}, the coupling is handled 
explicitly, that is, first we solve the problem featuring the organic carbon with 
fixed temperature then we approximate the heat equation using the information 
arising from the previously computed $\Corg$. \\
Within this framework, at time $t=t_n$ the conditions of humidity and temperature 
for the organic carbon equation read as
$$
\Psi_1(w) \equiv \frac{\wMax}{4}
,\qquad
\Psi_2(T) = \max \left(0, 1-\frac{|T_n-\Topt|}{A_T} \right)
$$
where $T_n$ is the temperature at the previous iteration.

\begin{figure}[hbtp]
\centering
\subfloat[Evolution of the organic carbon.]
{
\includegraphics[width=0.48\columnwidth]{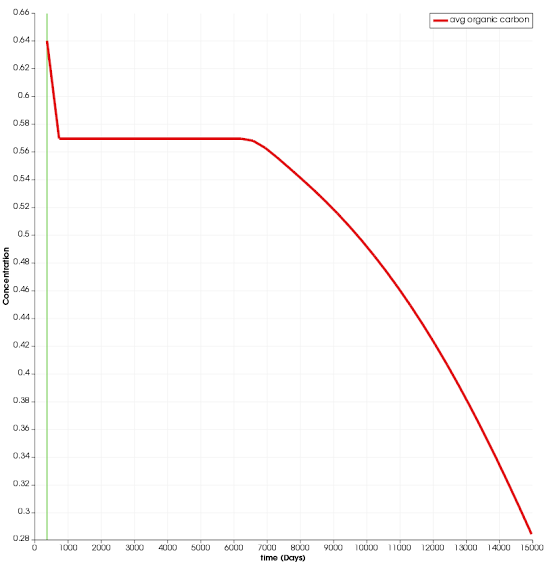}
\label{fig:coupled_carb_ev}
}
\hfill
\subfloat[Evolution of the temperature.]
{
\includegraphics[width=0.48\columnwidth]{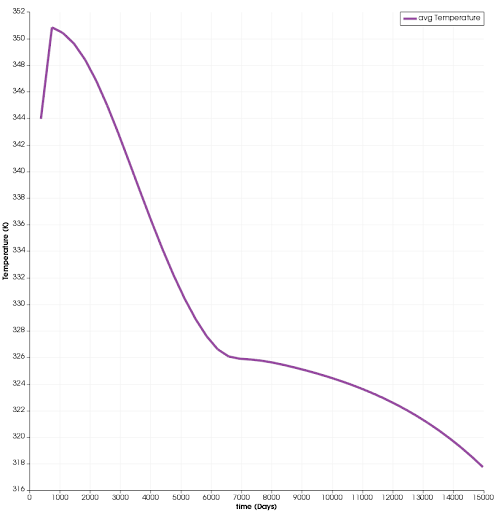}
\label{fig:coupled_temp_ev}
}
\caption{Evolution of the quantity of organic carbon and the temperature inside the alveolus 
between $t=0$ and $t=40$ years.}
\label{fig:coupled_temp_carb_ev}
\end{figure}

As in the previous case, we consider an initial concentration of organic carbon equal to 
$1$ and we observe it decreasing in figure \ref{fig:coupled_carb_ev} due to the 
bacterial activity. 
We verify that the quantity of organic carbon inside the bioreactor landfill 
decays towards zero during the lifetime of the facility. At the 
same time, the temperature increases as a result of the methanogenic process catalyzed 
by the microbiota (Fig. \ref{fig:coupled_temp_ev}). 
Nevertheless, when the temperature goes beyond the tolerated variation $A_T$, the 
second condition in (\ref{eq:psi1psi2}) is no more fulfilled and the chemical 
reaction is prevented. 
We may observe this behavior in figures \ref{fig:coupled_carb_ev}-\ref{fig:coupled_temp_ev} 
between $t=3$ years and $t=20$ years. 
Once the temperature is inside the admissible range $[\Topt - A_T; \Topt + A_T]$ 
again (starting approximately from $t=20$ years), the reaction (\ref{eq:chemical_reaction}) is allowed, 
the organic carbon is consumed and influences the temperature which slightly increases 
again before eventually decreasing towards the end-of-life of the bioreactor.
\begin{figure}[p]
\centering
\vspace{-20pt}
\subfloat[Lifetime: $1$ year.]
{
\includegraphics[width=0.5\columnwidth]{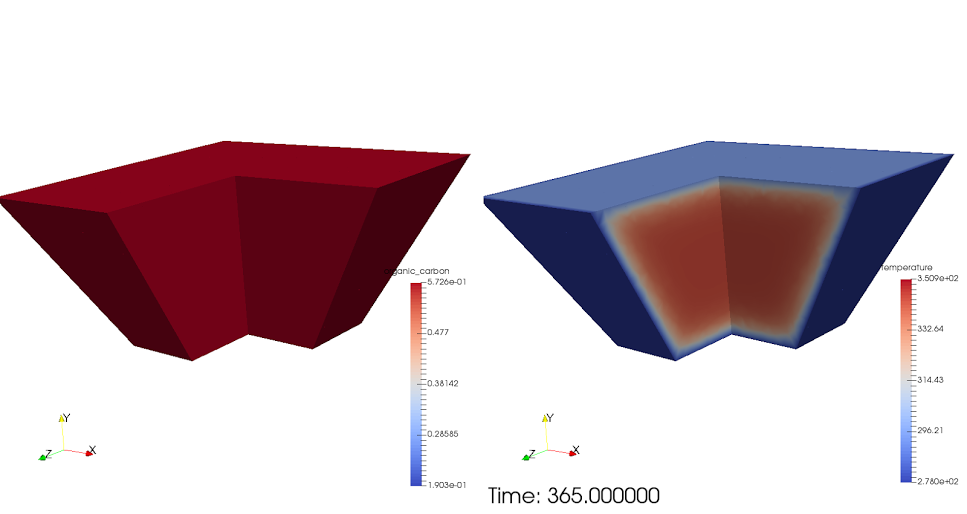}
\label{fig:coupled1}
}
\vspace{-10pt}\\
\subfloat[Lifetime: $10$ years.]
{
\includegraphics[width=0.5\columnwidth]{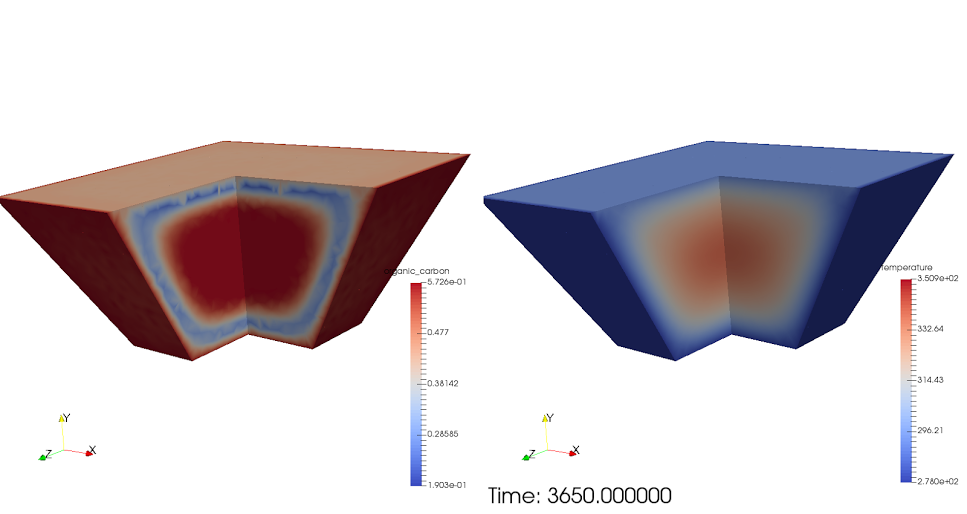}
\label{fig:coupled10}
}
\vspace{-10pt}\\
\subfloat[Lifetime: $20$ years.]
{
\includegraphics[width=0.5\columnwidth]{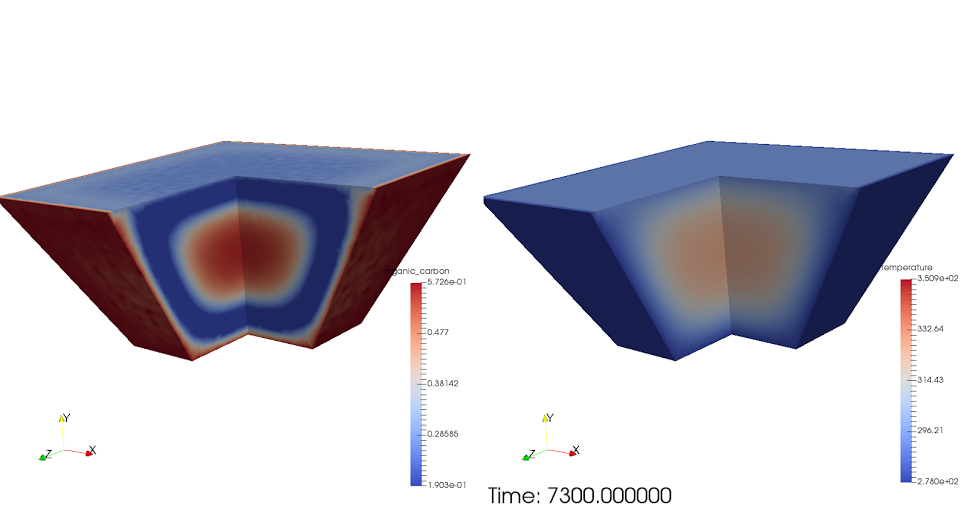}
\label{fig:coupled20}
}
\vspace{-10pt}\\
\subfloat[Lifetime: $40$ years.]
{
\includegraphics[width=0.5\columnwidth]{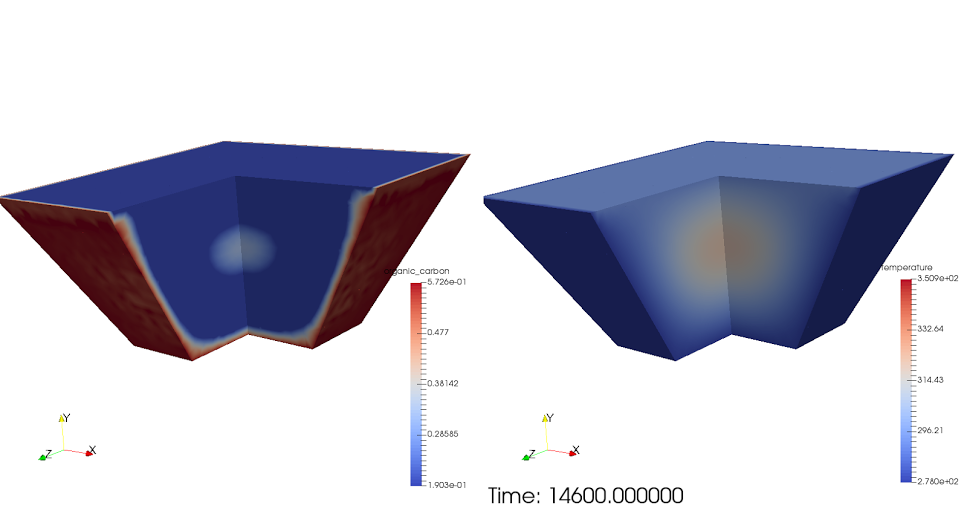}
\label{fig:coupled40}
}
\caption{Coupled evolution of the organic carbon (left) and the temperature (right) 
inside the alveolus at $t=1$ year, $t=10$ years, $t=20$ years and $t=40$ years.}
\label{fig:coupled_snapshots}
\end{figure}
Eventually, in figure \ref{fig:coupled_snapshots} we report some snapshots of the 
solutions of the coupled system (\ref{eq:organic_carbon_eq})-(\ref{eq:heat_eq}).

\section{Conclusion}
\label{ref:conclusion}

In this work, we proposed a first attempt to mathematically model the physical
and chemical phenomena taking place inside a bioreactor landfill. A set of $7$
coupled equations has been derived and a Finite Element discretization has been
introduced using \Feel.
A key aspect of the discussed model is the tuning of the coefficients appearing
in the equations. On the one hand, part of these unknowns represents physical
quantities whose values may be derived from experimental studies. On the other
hand, some parameters are scalar factors that have to be estimated by means of
heuristic approaches. A rigorous tuning of these quantities represents a
major line of investigation to finalize the implementation of the model in the 
SiViBiR++ module and its validation with real data.

The present work represents a starting point for the development of
mathematically-sound investigations on bioreactor landfills.
From a modeling point of view, some assumptions may be relaxed, for
example by adding a term to account for the death rate of bacteria or the
cooling effect due to water injection inside the bioreactor.
The final goal of SiViBiR++ project is the simulation of long-time behavior
of
the bioreactor in order to perform forecasts on the methane production and
optimize the control strategy. The associated inverse problems and
PDE-constrained optimization problems are likely to be numerically
intractable due to their complexity and their dimension thus the study of
reduced order models may be necessary to decrease the overall computational
cost.

\appendix
\section{Summary of the unknown parameters}
\label{sec:parameters}

In the following table, we summarize the values of some unknown parameters which 
were deduced during the present work. We highlight that all these quantities have 
been estimated via heuristic approaches and a rigorous verification/validation 
procedure remains necessary before their application to real-world problems.

\setlength\LTleft{0pt}
\setlength\LTright{0pt}

\begin{longtable}{@{\extracolsep{\fill}}rlll}
Parameter & Description & Value & Unit \\
\hline
\hline
\\
\endhead

$\Phi$           & Porosity of the medium                      & \SI{0.3}{}     & \\
$D$              & Permeability                                & \SI{e-11}{}    & \SI{}{m^2}  \\
$b_0$            & Initial concentration of bacteria           & \SI{1.0}{}     & \\
$\Corg_0$        & Initial concentration of organic carbon     & \SI{1.0}{}     & \\
$a_b$            & Rate of consumption of organic carbon       & \SI{e-5}{}     & \SI{}{m^6\kg^{-2}\day^{-1}} \\
$c_b$            & Rate of creation of bacteria                & \SI{1.0}{}     & \\
\\
\hline
\\
$\Topt$          & Optimal temperature for the reaction        & \SI{308}{}     & \SI{}{K} \\
$A_T$            & Tolerated variation of temperature          & \SI{20}{}      & \SI{}{K} \\
$c_T$            & Rate of production of heat by the chemical reaction          & \SI{e2}{}      & \SI{}{K} \\
$k_T$            & Thermal conductivity                        & \SI{9e-2}{}    & \SI{}{m^{2}\day^{-1}} \\
$T_g$            & Temperature of the soil                     & \SI{278}{}     & \SI{}{K} \\
$T_m$            & Temperature of the geomembrane              & \SI{293}{}     & \SI{}{K} \\
$T_0$            & Initial temperature                         & \SI{293}{}     & \SI{}{K} \\
\\
\hline
\\
$\wMax$          & Maximal admissible quantity of water        & \SI{100}{}     & \SI{}{\kg.m^{-3}} \\
$\|u_w\|$        & Velocity of the water                       & \SI{2.1}{}     & \SI{}{m.\day^{-1}} \\
$k_w$            & Diffusion coefficient of the water          & \SI{8.6e-2}{}  & \SI{}{m^2\day^{-1}} \\
$w_0$            & Initial quantity of water                   & \SI{50}{}      & \SI{}{\kg.m^{-3}} \\
\\
\hline
\\
$H_0$            & Constant for the vapor pressure            & \SI{-9.56e4}{} & \SI{}{\Pa} \\
$H_1$            & Constant for the vapor pressure            & \SI{337.89}{}  & \SI{}{\Pa.K^{-1}} \\
$c_{h\to w}$     & Condensation rate                           & \SI{e-1}{}     & \SI{}{\day^{-1}} \\
$\mu_\text{gas}$ & Dynamical viscosity of the gas              & \SI{1,3}{}     & \SI{}{Pa.\day^{-1}} \\
$c_M$            & Rate of production of methane               & \SI{1.8e7}{}   & \SI{}{\kg.m^{-3}} \\
$c_C$            & Rate of production of carbon dioxide        & \SI{2.6e7}{}   & \SI{}{\kg.m^{-3}} \\
$c_h$            & Rate of production of water vapor           & \SI{2.5e6}{}   & \SI{}{\kg.m^{-3}} \\
$M_0$            & Initial concentration of methane           & \SI{1.0}{}     & \\
$\Cdx_0$         & Initial concentration of carbon dioxide     & \SI{1.0}{}     & \\
$h_0$            & Initial concentration of water vapor        & \SI{1.0}{}     & \\
\\
\hline
\\
$\Jout$          & Outgoing flow of biogas                     & \SI{258}{}     & \SI{}{m^3\day^{-1}} \\
$\Jin$           & Incoming flow of water and leachates         & \SI{258}{}     & \SI{}{m^3\day^{-1}} \\
\\
\hline
\hline
\\
\caption{Summary of the parameters involved in the $7$-equations model.}
\end{longtable}

\subsubsection*{Acknowledgments}

The authors are grateful to Alexandre Ancel (Universit\'e de Strasbourg) and the
\Feel community for the technical support.
The authors wish to thank the CEMRACS 2015 and its organizers.

\bibliographystyle{abbrv}
\bibliography{./sivibirBiblio}
\end{document}